\documentclass[11pt,a4paper]{article}
\usepackage{amsmath, amssymb}
\usepackage{latexsym, color}
\usepackage{epsfig}
\usepackage{ascmac}
\usepackage{amscd} 
\usepackage{comment}
\usepackage{hyperref}
\hypersetup{
    colorlinks=true,
    linkcolor=blue,
    citecolor=blue,
    urlcolor=blue
}
\numberwithin{equation}{section}

\newtheorem{thm}{Theorem}[section]
\newtheorem{cor}[thm]{Corollary}
\newtheorem{lem}[thm]{Lemma}
\newtheorem{prop}[thm]{Proposition}
\newtheorem{definition}[thm]{Definition}

\newenvironment{pf}[1][Proof]
{\par\noindent\textit{#1}. \ignorespaces}
{\hfill\ensuremath{\square}\par}
\newcommand{\p}{\partial}

\newcommand{\Cbb}{\mathbb{C}}

\newcommand{\Nbb}{\mathbb{N}}

\newcommand{\Rbb}{\mathbb{R}}
\newcommand{\Sbb}{\mathbb{S}}
\newcommand{\Zbb}{\mathbb{Z}}

\newcommand{\la}{\langle}
\newcommand{\ra}{\rangle}

\newcommand{\Hcal}{\mathcal{H}}

\newcommand{\Kcal}{\mathcal{K}}

\newcommand{\Scal}{\mathcal{S}}

\def\Bx{{\bf x}}
\def\By{{\bf y}}

\newcommand{\Ga}{\alpha}

\newcommand{\Gvf}{\varphi}

\newcommand{\Gl}{\lambda}

\newcommand{\Gm}{\mu}
\newcommand{\Gv}{\nu}

\newcommand{\Gt}{\theta}

\newcommand{\Gr}{\rho}

\newcommand{\Gx}{\xi}
\newcommand{\Gy}{\psi}

\newcommand{\GG}{\Gamma}

\newcommand{\GO}{\Omega}

\DeclareMathOperator{\sub}{sub}

\newcommand{\beq}{\begin{equation}}
\newcommand{\eeq}{\end{equation}}

\begin{document}
\title{Precise Eigenvalue Asymptotics \\ for Electrostatic Layer Potential Operators \thanks{\footnotesize Y.-G. Ji was supported by KIAS Individual Grant No. HP089001 via the June E Huh Center for Mathematical Challenges at Korea Institute for Advanced Study. Y. Miyanishi was supported by JSPS (of Japan) KAKENHI Grant Numbers JP25K07032, JP25K07081.}}

\author{
Yong-Gwan Ji\thanks{June E Huh Center for Mathematical Challenges, Korea Institute for Advanced Study, Seoul 02455, S. Korea (ygji@kias.re.kr.)} 
\and Daisuke Kawagoe\thanks{Graduate School of Informatics, Kyoto University, Kyoto 606-8501, Japan (d.kawagoe@acs.i.kyoto-u.ac.jp).}  
\and Yoshihisa Miyanishi\thanks{Department of Mathematical Sciences, Faculty of Science, Shinshu University, Matsumoto 390-8621, Japan (miyanishi@shinshu-u.ac.jp).} }

\maketitle

\begin{abstract}
We investigate the eigenvalue asymptotics of electrostatic layer potential operators associated with the Laplace equation on smooth bounded domains in two and three dimensions. In two dimensions, we establish precise eigenvalue asymptotics for the single-layer potential operator, up to a remainder of order $O(j^{-\infty})$. In three dimensions, we obtain refined eigenvalue asymptotics for the single-layer potential operator on arbitrary smooth surfaces and for the Neumann--Poincar\'e operator on smooth surfaces 
with everywhere positive Gaussian curvature. Under suitable nonperiodicity assumptions on the associated Hamiltonian flows, the remainder estimates in three dimensions improve to $o(j^{-1})$. 
\end{abstract}


\noindent{\footnotesize {\bf Key words}. Layer potentials, single-layer potential operator, Neumann--Poincar\'e operator, \\ eigenvalue asymptotics, Weyl's law}
\renewcommand{\contentsname}{Table of Contents }
\tableofcontents

\section{Introduction}
Let $\Omega$ be a bounded open set with Lipschitz boundary $\Gamma=\p\GO$ in $\Rbb^d\ (d=2, 3)$. We take the fundamental solution of the Laplace operator $\Delta$ to be 
\beq\label{eq: fundamental sol of Laplacian}
E(x, y)=
\begin{cases}
\frac{1}{2\pi} \log |x-y|\quad d=2, \\
-\frac{1}{4\pi} \frac{1}{|x-y|}\quad\quad\hspace{3mm} d=3.
\end{cases}
\eeq
Our main purpose is to investigate the precise asymptotic behavior 
of the eigenvalues of two boundary integral operators associated with layer potentials.  

The first boundary integral operator is the single-layer potential operator associated with $\Gamma=\p\GO$, 
which is defined by 
\beq
\Scal_{\Gamma} [f](x) := \int_{\Gamma} E(x, y)f(y)\; dS_{y}, \quad x \in \GG,
\eeq
where $f \in L^2 (\Gamma)$ and $dS_y$ is the line or surface element. 
Since the single-layer potential operator defines a bounded map $\mathcal S_{\Gamma}:H^{-1/2}(\Gamma)\rightarrow H^{1/2}(\Gamma)$ for Lipschitz boundaries and the embedding \(H^{1/2}(\Gamma)\hookrightarrow L^{2}(\Gamma)\) is compact, the operator
$\mathcal S_{\Gamma}: L^2(\Gamma)\rightarrow L^2(\Gamma)$ is compact. 
Moreover, since its integral kernel is real-valued and symmetric, $\Scal_\Gamma$ is self-adjoint on $L^2(\Gamma)$ (see, e.g., \cite{C,Z}).
We emphasize that the spectral properties of
single-layer potential operators associated with
the Laplace and Helmholtz equations play an
important role in various applications, including
electrostatic capacity \cite{BAD}, numerical methods
for boundary integral equations \cite{G}, and acoustic scattering \cite{AT, TA}.

The second important boundary integral operator is the Neumann--Poincar\'e (abbreviated by NP) operator, which is defined by   
\beq\label{NP operator}
\Kcal_{\Gamma} [f](x)=\mbox{p.v.}\int_{\Gamma} \p_{\nu_y} E(x, y)f(y)\; dS_{y}, \quad x \in \GG,
\eeq
where $\p_{\nu_y}$ is the outer normal derivative at $y\in \Gamma$. Here the notation ``p.v.'' denotes the principal value. The NP operator (also called the double-layer potential operator) is bounded when $\Gamma$ is Lipschitz and compact when $\Gamma$ is of class $C^{1,\alpha}$ for some $\alpha>0$ (see, e.g., \cite{AKMP, FKM, MR2308861}).
The spectral properties of the NP operator also play an important role in various applications, including plasmonic resonance, cloaking, stress concentration, and inverse problems (see, e.g., \cite{ACKLM, ACKLY, ADM, AK, AKMN, AKMP2}).

We note here that different sign conventions are used in the literature for the fundamental solution of the Laplacian and, consequently, for the integral kernel of the NP operator. Accordingly, statements concerning positive and negative eigenvalues may be interchanged depending on the convention adopted. 
We adopt the sign convention used in \cite{FKM}, with the fundamental solution given by \eqref{eq: fundamental sol of Laplacian}. We also define the second fundamental form using the outward unit normal, so that the
principal curvatures of a sphere are negative under our convention.

 {\it In what follows, we assume that all boundaries are  smooth}. Although this assumption may be somewhat restrictive, we impose it for technical reasons in order to study the precise asymptotic behavior of the eigenvalues. For clarity, we fix the following convention for ordering eigenvalues. For a compact operator $A$ with real eigenvalues, we denote its nonzero eigenvalues, arranged in nonincreasing order of their absolute values and counted with multiplicity, by $\lambda_j(A)$. For a positive self-adjoint elliptic operator $A$ with compact resolvent, we use the same notation $\lambda_j(A)$ for its eigenvalues, arranged in nondecreasing order and counted with multiplicity. 
Under this smoothness assumption, we obtain refined eigenvalue asymptotics for the single-layer potential operator and the NP operator. We treat the two- and three-dimensional cases separately. 

To state our results in two dimensions, we first recall that the kernel of the single-layer potential operator $\Scal_\Gamma$  has dimension at most one. Indeed, if $\varphi\in\ker \Scal_\Gamma$ satisfies $\int_\Gamma\varphi\, dS=0$, then the corresponding single-layer potential tends to zero at infinity. The uniqueness of the interior and exterior Dirichlet problem 
with decay at infinity, together with the jump relation, implies that $\varphi=0$. Consequently, 
the linear functional \(\varphi\mapsto\int_\Gamma\varphi\,dS\) is injective on $\ker \Scal_\Gamma$. Since its range is contained in a one-dimensional space, we conclude that $\dim\ker \Scal_\Gamma \leq1.$ Moreover, it is known that $S_\Gamma$ has a nontrivial kernel if and only if the logarithmic capacity of $\Gamma$ is equal to one (see, e.g., \cite{DH, Mc, S}), that is, 
$$
\dim\ker \Scal_\Gamma=
\begin{cases}
1,&\operatorname{cap}(\Gamma)=1,\\
0,&\operatorname{cap}(\Gamma)\neq1.
\end{cases}
$$
Then our main results are stated in the following three theorems. We begin with the two-dimensional case. 
\begin{thm}\label{thm: single-layer accurate 2D eigenvalues}
Let $\Omega \subset \Rbb^2$ be a bounded open set with smooth boundary $\Gamma=\p\GO$. Let $\Gamma_1, \cdots, \Gamma_k$ be the connected components of the boundary $\Gamma=\p\GO$, with lengths $l(\Gamma_i), 1\leq i \leq k$. Let $R=\left( \frac{1}{l(\Gamma_1)}, \cdots, \frac{1}{l(\Gamma_k)} \right)$ be a finite sequence of positive numbers, where repetitions are allowed. Then the nonzero eigenvalues of the single-layer
potential operator $\Scal_\Gamma$, arranged in nonincreasing order of their absolute values and
counted with multiplicity, satisfy 
\beq\label{eq: S(R) ordered eigenvalues of single layer potential}
\lambda_{j}(\Scal_{\Gamma})=-\frac{1}{4\pi}({S(R)}_{j+\dim\ker \Scal_{\Gamma}})^{-1}+O(j^{-\infty}) \quad \mbox{as}\ j\rightarrow \infty.
\eeq
Here $S(R)_j$ denotes the $j$-th smallest element of the multiset $S(R)$, counting multiplicities
and including zeros, as defined in Section \ref{sec: DN map eigenvalues notation}. 

In particular, if $\Gamma$ is connected and has length $l(\Gamma)$, then
\beq\label{eq: single layer dimker Scal=0 case}
\lambda_{2j}(\Scal_\Gamma)=\lambda_{2j+1}(\Scal_\Gamma)+O(j^{-\infty}) = -\frac{l(\Gamma)}{4\pi j} +O(j^{-\infty}) \quad \mbox{as}\ j\rightarrow \infty
\eeq
when $\ker \Scal_\Gamma=\{0\}$. If $\dim\ker S_\Gamma=1$, then
\beq\label{eq: single layer dimker Scal=1 case} 
\lambda_{2j-1}(\Scal_\Gamma) = \lambda_{2j}(\Scal_\Gamma) +O(j^{-\infty}) =
-\frac{l(\Gamma)}{4\pi j} +O(j^{-\infty}) \quad \mbox{as}\ j\rightarrow \infty.
\eeq
\end{thm}
For the two-dimensional NP operator, although such precise asymptotic behavior is not yet known, several aspects of its eigenvalue behavior are studied in recent work involving the first and third authors (see, e.g., \cite{AKMY, FKM}.). 

In the three-dimensional case, we obtain the following results. 
\begin{thm}[Weyl's law for the single-layer potential operator]\label{thm: Weyl asymptotic of single-layer potentials in three dimensions}
Let $\Omega \subset \Rbb^3$ be a smooth bounded open set with smooth boundary $\Gamma$. Then 
\beq\label{eq: single-layer rough spectral asymptotic of eigenvalues in 3 dimension}
\lambda_{j}(\Scal_{\Gamma})=-\frac{1}{4\sqrt{\pi}} {\text{\rm Area}(\Gamma)}^{1/2}  j^{-1/2}+O(j^{-1}) \quad \mbox{as}\ j\rightarrow \infty.
\eeq
Furthermore, if the set of periodic points of the geodesic flow on the cosphere bundle $S^*\Gamma$ has Liouville measure zero, then 
\beq\label{eq: second term of single-layer}
\lambda_{j}(\Scal_{\Gamma})=-\frac{1}{4\sqrt{\pi}} {\text{\rm Area}(\Gamma)}^{1/2}  j^{-1/2}+o(j^{-1}) \quad \mbox{as}\ j\rightarrow \infty.
\eeq
\end{thm}
It is known that, for generic smooth surfaces,
the set of periodic points of the geodesic flow
on the cosphere bundle has Liouville measure zero
(see, e.g., \cite{Cl, ST}).
Consequently, the remainder term in the eigenvalue
asymptotics of the single-layer potential operator
improves to $o(j^{-1})$ for generic smooth surfaces. 
In particular, there is no nonzero second-order term of order $j^{-1}$ in the eigenvalue asymptotics. 

We note that the conclusion of Theorem \ref{thm: Weyl asymptotic of single-layer potentials in three dimensions} was previously announced by the third author in \cite{AKMP}. The purpose of the present paper is to provide detailed proofs of \eqref{eq: single-layer rough spectral asymptotic of eigenvalues in 3 dimension} and \eqref{eq: second term of single-layer}, based primarily on microlocal analysis and the Dirichlet-to-Neumann (abbreviated by DtN) map, together with explicit examples.

For the NP operator, detailed results in three dimensions can be obtained 
for the smooth surfaces with everywhere positive Gaussian curvature. 
Since the Gaussian curvature of $\Gamma$ is everywhere positive, the mean curvature is nowhere zero 
and has a constant sign on each connected component of $\Gamma$.
Thus we can decompose $\Gamma$ as
\beq
\Gamma=\Gamma^+\sqcup\Gamma^-,
\eeq
where $\Gamma^+$ and $\Gamma^-$ are the unions of the connected components on which
the mean curvature is positive and negative, respectively.
\begin{thm}[Weyl's law for the NP operator]\label{thm: Weyl's law for NP}
Let $\Omega \subset \Rbb^3$ be a bounded open set with smooth boundary $\Gamma$, and assume that
the Gaussian curvature of $\Gamma$ is everywhere positive. Then the components of $\Gamma$ can be
classified according to the sign of their mean curvature. We write
\beq
\Gamma = \Gamma^{+} \sqcup \Gamma^{-} 
=\left(\bigsqcup_{k=1}^{m}\Gamma_k^+\right)
\sqcup
\left(\bigsqcup_{l=1}^{n}\Gamma_l^-\right),
\eeq
where $\Gamma_k^+$ and $\Gamma_l^-$ denote
the connected components on which the mean curvature
is everywhere positive and everywhere negative,
respectively.
We suppose that the set of periodic points of the Hamiltonian flow generated by $p=|K_0|^{-1}$ on the energy surface
\beq
\{(x,\xi)\in T^*\Gamma\setminus 0 \mid p(x,\xi)=1\}
\eeq
has Liouville measure zero. Here the function $K_0$ on $T^{*}\Gamma$ is defined in Section \ref{sec: symbol}.  We denote $\lambda_j^{\pm}(\Kcal_{\Gamma})$ as the positive and
negative eigenvalues of $\Kcal_{\Gamma}$, respectively, arranged in nonincreasing order of their absolute
values and counted with multiplicity. Then we have 
\beq\label{eq: accurete positive eigenvalue asymptotics}
\begin{aligned}
\lambda_j^{+}(\Kcal_\Gamma) &= \left( \frac{3W(\Gamma^-)-2\pi\chi(\Gamma^-)} {128\pi} \right)^{1/2} j^{-1/2} +o(j^{-1}) \\
&=\left\{\sum_{l=1}^{n} \left( \frac{3W(\Gamma_l^-)-2\pi\chi(\Gamma_l^-)} {128\pi} \right) \right\}^{1/2} j^{-1/2} +o(j^{-1}) 
\quad \text{as}\  j\to\infty,
\end{aligned}
\eeq

Moreover if $\Gamma^{+} \neq \varnothing$, 
\beq
\begin{aligned}
\lambda_j^{-}(\Kcal_\Gamma) &= -\left( \frac{3W(\Gamma^+)-2\pi\chi(\Gamma^+)} {128\pi} \right)^{1/2} j^{-1/2} +o(j^{-1}) \\
&= -\left\{\sum_{k=1}^{m} \left( \frac{3W(\Gamma_k^+)-2\pi\chi(\Gamma_k^+)} {128\pi}\right) \right\}^{1/2} j^{-1/2} +o(j^{-1})
\quad \text{as}\  j\to\infty. 
\end{aligned}
\eeq
Here the notation $W(\Gamma)=\int_\Gamma H^2\,dS$ denotes the Willmore energy, where $H=\frac{\kappa_1+\kappa_2}{2}$ is the mean curvature of $\Gamma$, and $\chi(\Gamma)$ denotes its Euler characteristic.  The Neumann--Poincar\'e (NP)  flow generated by $p=|K_0|^{-1}$ and its periodic points are defined in Definitions \ref{def: NP Hamilton flow} and \ref{def: NP periodic points}, respectively. 
\end{thm}
We emphasize that the NP operator has only finitely many negative eigenvalues when $\Gamma^{+}=\varnothing$ (see Section \ref{sec: Eigenvalues of the NP operator in three dimensions}). Thus, in this case, we focus on the asymptotic behavior of the positive eigenvalues  \eqref{eq: accurete positive eigenvalue asymptotics}. 

We also note that, for general $C^{2,\alpha}$ ($\alpha>0$) surfaces without any positivity assumption on the Gaussian curvature, the known Weyl asymptotics yield a remainder estimate of order $o(j^{-1/2})$ (see, e.g., \cite{Miyanishi:Weyl, MR}). For smooth surfaces with everywhere positive Gaussian curvature, we improve the remainder estimate to $o(j^{-1})$, under the additional assumption that the set of periodic points of the associated Hamiltonian flow (the NP flow) has Liouville measure zero.

To establish these main theorems, the paper is organized as follows. In Section \ref{sec: prelminaries}, we review the basic properties of layer potential operators and their relations with the Dirichlet-to-Neumann (DtN) map. We also discuss layer potential operators from the viewpoint of pseudodifferential operator theory, including their principal and subprincipal symbols. In addition, we introduce the Neumann--Poincar\'e (NP) flow as the Hamiltonian flow generated by the reciprocal of the absolute value of the principal symbol of the NP operator. We then introduce the periodic points of this flow on the unit energy surface.

Based on the preliminary results established in Section \ref{sec: prelminaries}, we prove in Section \ref{sec: Eigenvalue asymptotics of layer potential operators for C^infty} our main results, Theorems \ref{thm: single-layer accurate 2D eigenvalues}, \ref{thm: Weyl asymptotic of single-layer potentials in three dimensions}, and \ref{thm: Weyl's law for NP}, which give refined descriptions of the eigenvalue asymptotics. The proofs are presented separately for the two- and three-dimensional cases.

\section{Preliminaries}\label{sec: prelminaries}

\subsection{Jump relation and Plemelj's principle}\label{sec: Jump relation and Plemelj}

For completeness, we review the basic properties in this section. We first define the single-layer potential operator $\Scal_{\Omega}$, viewed as a mapping from functions on the boundary to functions in $\Rbb^d \setminus \Gamma$, as follows.
\beq
\Scal_{\GO}[f](x) = \int_{\Gamma} E(x,y)f(y)\,dS_y,
\qquad
x\in\mathbb{R}^d\setminus\Gamma.
\eeq
Thus the single-layer boundary operator $\mathcal{S}_{\Gamma}$ considered in this paper is then obtained as the Dirichlet trace of $\mathcal{S}_{\Omega}$. In other words, $\mathcal{S}_{\Gamma}$ is the restriction of $\mathcal{S}_{\Omega}$ to the boundary $\Gamma$ in the sense of traces.

Then one of the fundamental relations is the well-known jump formula:     
\beq\label{gs-smooth}
\p_\Gv^{-} \Scal_{\GO}[\phi]|_{\Gamma}=\left(-\frac{1}{2} + \Kcal^{*}_\Gamma \right)[\phi].
\eeq
Here, $\partial_\nu^-$ denotes the interior non-tangential limit of the normal derivative (see, e.g., \cite{Kr, Mc}).
It is also well known that Plemelj's symmetrization principle, also referred to as Calder\'on's identity, states that
\beq
\Scal_{\Gamma} \Kcal_{\Gamma}^* = \Kcal_{\Gamma} \Scal_{\Gamma},
\eeq
where $\Kcal_{\Gamma}$ denotes the $L^2(\Gamma)$-adjoint of $\Kcal_{\Gamma}^*$. For further details, see \cite{MR2308861}.

We denote by $H^s = H^s(\Gamma)$ and $\| \cdot \|_s$, $s \in \Rbb$, the usual Sobolev space on $\Gamma$ and its norm, respectively; in particular, $H^0(\Gamma) = L^2(\Gamma)$.
We define a bilinear form $\langle \cdot, \cdot \rangle_{\Hcal^{*}}$ on $H_0^{-1/2}(\Gamma) := \{ \Gvf \in H^{- 1/2}(\Gamma); \la \Gvf, 1 \ra_{H^{-1/2}, H^{1/2}} = 0 \}$ by
\beq
  \la \Gvf, \Gy \ra_{\Hcal^*} := - \la \Gvf, \Scal_{\Gamma} [\Gy] \ra_{H^{-1/2}, H^{1/2}} \label{eq:new_inner_product1}
\eeq
where the pairing on the right-hand side is
understood in the sense of Sobolev duality
between $H^{-1/2}(\Gamma)$ and $H^{1/2}(\Gamma)$.
We note that the right hand side of \eqref{eq:new_inner_product1} is well-defined since $\Scal_{\Gamma}$ maps $H^{- 1/2}$ to $H^{1/2}$.
In fact, the well-known positivity of $-S_\Gamma$ on the mean-zero subspace implies that $\la \cdot, \cdot \ra_{\Hcal^*}$ defines an inner product on $H_0^{-1/2}(\Gamma)$ (see, e.g., \cite{Kang}).
We put $\Hcal_0^* := H_0^{-1/2}$ equipped with the inner product $\la \cdot, \cdot \ra_{\Hcal^*}$.
We can easily see that $\Kcal_{\Gamma}^*$ can be symmetrized on $\Hcal_0^*$ by Plemelj's symmetrization principle, {\it i.e.}
\beq
  \la \varphi, \Kcal_{\Gamma}^*[\psi] \ra_{\Hcal^*} = \la \Kcal_{\Gamma}^*[\varphi], \psi \ra_{\Hcal^*}. \label{eq:Plemelj_0}
\eeq

Furthermore, it is proved in \cite{MR0222317} that the eigenvalues of $\Kcal_{\Gamma}^*$ on $L^2(\Gamma)$ lie in $[-1/2, 1/2]$.
It is easy to see that $\Kcal_{\Gamma}^*$ is compact on $L^2$, if $\Gamma$ is smooth.
Moreover, we can see in \cite{MR2308861} that $\Kcal_{\Gamma}^*$ is also compact on $H^{-1/2}$ and has the same eigenvalues and eigenfunctions as on $L^2$.
We suppose that $\Omega$ is a bounded open set with
smooth boundary $\Gamma$, consisting of $n$
connected components $\Omega_1,\ldots,\Omega_n$.
We write $\Gamma_j=\partial\Omega_j$ and
denote by $\chi_j$ the function on $\Gamma$
which is equal to $1$ on $\Gamma_j$ and
$0$ on the remaining components.

It is known that, if $\varphi$ is an eigenfunction of $\Kcal_\Gamma^*$ corresponding to the eigenvalue $1/2$, then
$\Scal_\Omega[\varphi]$ is constant on each connected component of $\Omega$.
In two dimensions, these constants may all vanish even when $\varphi\neq 0$.
Consequently, the pairing induced by $-\Scal_\Gamma$ may fail to be positive definite on $H^{-1/2}(\Gamma)$.

We here introduce the closed subspace
\beq
H_{j,0}^{-1/2}(\Gamma) = \left\{ \varphi\in H^{-1/2}(\Gamma)|\ \la\varphi,\chi_j\ra=0 \right\} \quad (j=1,\ldots,n). 
\eeq
Since the $n$ conditions are linearly independent, we have
\beq
\dim\bigl(
H^{-1/2}(\Gamma)/
\bigcap_{j=1, \ldots, n} H_{j,0}^{-1/2}(\Gamma)
\bigr)=n.
\eeq
It follows from $\Kcal_\Gamma[\chi_j]
=\frac12\chi_j \ (j=1,\ldots,n) $ and the Fredholm alternative implies that
\beq
\dim\ker(\Kcal_\Gamma^*-\tfrac12 I) \geq n.
\eeq

For $ \displaystyle\varphi\in \ker(\Kcal_\Gamma^*-\tfrac12 I)
\cap \bigcap_{j=1, \ldots, n} H_{j,0}^{-1/2}(\Gamma)$, the jump relation implies that
$\Scal_\Omega[\varphi]$ has vanishing interior normal derivative on $\Gamma$.
Hence, $\Scal_\Omega[\varphi]$ is constant on each connected component $\Omega_j$.
It follows from continuity of the single-layer potential that there exist constants $c_j$ such that
\beq
\Scal_\Gamma[\varphi] = \sum_{j=1}^n c_j\chi_j.
\eeq
Since
$\displaystyle\varphi\in \bigcap_{j=1, \ldots, n} H_{j,0}^{-1/2}(\Gamma)$, we obtain
\beq
\la\varphi,\Scal_\Gamma[\varphi]\ra = \sum_{j=1}^n c_j\la\varphi,\chi_j\ra =0.
\eeq
Then the positive definiteness of $-\Scal_\Gamma$
on $\bigcap_{j=1, \ldots, n} H_{j,0}^{-1/2}(\Gamma)$ yields $\varphi=0$.
Thus,
\beq
\ker(\Kcal_\Gamma^*-\tfrac12 I)
\cap \bigcap_{j=1, \ldots, n} H_{j,0}^{-1/2}(\Gamma) =\{0\}.
\eeq
Together with the preceding dimension estimate, this implies
\beq
\dim\ker(\Kcal_\Gamma^*-\tfrac12 I)=n.
\eeq
Consequently, the basis of $\ker(\Kcal_\Gamma^*-\tfrac12 I)$ is uniquely denoted as  
$\{\Gvf_1,\ldots,\Gvf_n\}$  satisfying 
\beq
\la\Gvf_i,\chi_j\ra=\delta_{ij}, \quad 1\leq i,j\leq n.
\eeq
Thus we have the direct-sum decomposition
\beq
H^{-1/2}(\Gamma) 
=\bigcap_{j=1, \ldots, n} H_{j,0}^{-1/2}(\Gamma) \oplus \operatorname{span} \{\Gvf_1,\ldots,\Gvf_n\}.
\eeq
Using this decomposition, we define a modified single-layer operator $\widetilde\Scal_\Gamma$ by 
\beq
\widetilde\Scal_\Gamma|_{\bigcap_{j=1, \ldots, n} H_{j,0}^{-1/2}(\Gamma)}= \Scal_\Gamma,
\quad
\widetilde\Scal_\Gamma[\Gvf_j] = -\chi_j, \quad j=1,\ldots,n, 
\eeq
and extend it linearly to $H^{-1/2}(\Gamma)$. We then define
\beq
\la\Gvf,\Gy\ra_{\Hcal^*}
=-\la\Gvf, \widetilde\Scal_\Gamma[\Gy] \ra_{H^{-1/2},H^{1/2}}
\eeq 
for $\Gvf,\Gy\in H^{-1/2}(\Gamma)$.

To verify positive definiteness of $-\widetilde\Scal_\Gamma$, we denote 
\beq
\Gvf=\psi+\sum_{j=1}^n a_j\Gvf_j, \quad \psi\in \bigcap_{j=1, \ldots, n} H_{j,0}^{-1/2}(\Gamma).
\eeq 
Since $\displaystyle\psi\in\bigcap_{j=1}^n H_{j,0}^{-1/2}(\Gamma)$, we have $\la\psi,\chi_j\ra=0. $ Moreover, since $\Scal_\Gamma[\Gvf_j]$ is constant on each $\Gamma_i$, the symmetry of $\Scal_\Gamma$ yields 
$\la\Gvf_j,\Scal_\Gamma[\psi]\ra = \la\Scal_\Gamma[\Gvf_j],\psi\ra =0.$ 
Thus, all cross terms vanish.
Consequently, 
\beq
\|\Gvf\|_{\Hcal^*}^2 =
-\la\psi,\Scal_\Gamma[\psi]\ra
+ \sum_{j=1}^n |a_j|^2.
\eeq
The positive definiteness of $-\Scal_\Gamma$ on $\displaystyle\bigcap_{j=1, \ldots, n} H_{j,0}^{-1/2}(\Gamma)$
shows that $\la\cdot, \cdot \ra_{\Hcal^*}$ defines an inner product on $H^{-1/2}(\Gamma)$.
Moreover, the standard coercivity estimate for $-\Scal_\Gamma$ on this subspace \cite{C, Steinbach},
together with the finite-dimensional direct-sum decomposition, yields 
\beq
\|\Gvf\|_{\Hcal^*} \approx
\|\Gvf\|_{H^{-1/2}(\Gamma)}.
\eeq 
We denote by $\Hcal^*$ the space $H^{-1/2}(\Gamma)$ equipped with this inner product.

Finally, the Plemelj symmetrization identity \eqref{eq:Plemelj_0} extends to
\beq
\widetilde\Scal_\Gamma\Kcal_\Gamma^* = \Kcal_\Gamma\widetilde\Scal_\Gamma.
\eeq
In fact, this identity holds on  $\displaystyle\bigcap_{j=1, \ldots, n} H_{j,0}^{-1/2}(\Gamma)$ by the original Plemelj identity 
and on $\operatorname{span}\{\Gvf_1,\ldots,\Gvf_n\}$ by the relations
$\Kcal_\Gamma^*[\Gvf_j] =\frac12\Gvf_j $ and $\Kcal_\Gamma[\chi_j]=\frac12\chi_j.$
Thus, $\Kcal_\Gamma^*$ is self-adjoint on $\Hcal^*$ (see, e.g., \cite{AK}).

In what follows, when $\Gamma$ is smooth, we occasionally use local coordinate charts to identify Sobolev spaces on $\Gamma$ 
with the corresponding Sobolev spaces on open subsets of $\mathbb{R}^{d-1}$. In particular, for surfaces in $\mathbb{R}^3$, 
we work in local coordinates on open subsets of $\mathbb{R}^2$.

\subsection{The Dirichlet--to--Neumann map and layer potentials}\label{sec: DN map and layer potentials}
When $\Scal_\Gamma : H^{-1/2}(\Gamma) \to H^{1/2}(\Gamma)$ is invertible, we can define $\varphi=\Scal_\Gamma^{-1}f$
for a given $f\in H^{1/2}(\Gamma)$.\ Then the harmonic extension of $f$ to $\Omega$ is represented by
\beq
u=\Scal_\Omega[\varphi].
\eeq
It follows from the interior jump relation that 
\beq
\Lambda f = \left(-\frac{1}{2} I+\Kcal_\Gamma^* \right)
\Scal_\Gamma^{-1}f.
\eeq
Here the notation $\Lambda$ denotes the Dirichlet--to--Neumann (abbreviated by DtN) map (see, e.g., \cite{GP, GPPS} for the details). 
Using either the self-adjointness of $\Lambda$ on $L^2(\Gamma)$
or Plemelj's symmetrization principle, we obtain the following Proposition. 
\begin{prop}\label{Relation between DN and layer potentials}
Let $\Omega$ be a bounded open set with Lipschitz boundary, and suppose that the single-layer operator $\Scal_\Gamma:H^{-1/2}(\Gamma) \to H^{1/2}(\Gamma)$ is invertible. Then the following identities hold.
\beq\label{eq: relation among single, DN and NP}
\left( -\frac{1}{2} I + \Kcal_{\Gamma}^* \right)\Scal^{-1}_{\Gamma}=\Scal^{-1}_{\Gamma}\left( -\frac{1}{2} I + \Kcal_{\Gamma} \right) = \Lambda.   
\eeq
\end{prop}
\subsection{Layer potentials as pseudo-differential operators in three dimensions}\label{sec: symbol}
In three dimensions,  $-\Scal_{\Gamma}$ is a positive, injective, self-adjoint elliptic pseudodifferential operator of order $-1$ on
$L^2(\Gamma)$. Moreover, it is an isomorphism from $H^{-1/2}(\Gamma)$ onto $H^{1/2}(\Gamma)$, and it satisfies Plemelj's symmetrization principle (also known as the Calder\'on's identity):
$$\Scal_{\Gamma} \Kcal^{*}_{\Gamma} =\Kcal_{\Gamma} \Scal_{\Gamma}.$$ Thus we can symmetrize the NP operator via single-layer potentials, namely,
\begin{align}\label{symmetrized NP}
  \widetilde{\Kcal_{\Gamma}} :={(-\Scal_{\Gamma})}^{-1/2} \Kcal_{\Gamma} {(-\Scal_{\Gamma})}^{1/2}
\end{align}
is  self-adjoint in $L^2(\Gamma)$. Our aim in this section is to find a pseudo-differential representation of the symmetrized NP operator $\widetilde{\Kcal_\Gamma}$ \eqref{symmetrized NP}. To this purpose, we represent the single-layer potential $\Scal_{\Gamma}$ and the NP operator ${\Kcal_{\Gamma}}$ as pseudo-differential operators (abbreviated by $\Gy$DOs).
Thanks to the smoothness condition, the asymptotic expansion of the integral kernel in  terms, homogeneous in $\Bx-\By$, yields a $\Gy$DO  for any local chart on the boundary manifold. The calculations in \cite{Miyanishi:Weyl} show that the principal symbol  of the NP operator, of order $-1$, equals in local co-ordinates $x',\xi$ 
\beq\label{eq: NP-PrincSymbol}
p_{\Kcal}(x',\xi)=-\frac{L(x') \xi_2^2 -2M(x')\xi_1 \xi_2 +N(x')\xi_1^2}{4 \det (g_{jk}) \left(\sum_{j, k}g^{jk}(x')\xi_j\xi_k\right)^{3/2}}
\eeq
(note especially here the minus sign in front of the fraction.) Here  $g_{jk}(x')$ denotes the metric tensor and $L(x'), M(x'), N(x')$ are the coefficients of the second fundamental form on the boundary $\Gamma$ in the local co-ordinates  $x'$ in some domain in $\Rbb^2$. It is important to emphasize that the unit normal vector used to define the second fundamental form is chosen to point
outward from $\GO$. See \cite{Ola} for the corresponding results in dimensions
three and higher.

It is well known \cite{AKMP, Taylor2000} that the single-layer potential operator is a classical $\Gy$DO of order $-1$ with principal symbol
$(-2 |\xi|_{x'})^{-1}.$ Here $|\xi|_{x'}=\sqrt{\sum_{j, k} g^{jk}(x')\xi_j \xi_k}$ on each local chart  (see, e.g., \cite{AKMP, Agranovich, MR}.).
We then find that its subprincipal symbol vanishes. 
In fact, for the geodesic normal coordinate centered at $X(0)=x_0 \in \Gamma$, it is known \cite{NK} that 
\beq\label{eq: distance in local chart geodesic normal coordinate}
|X(y) - X(0)| =|y|(1+O(|y|^2)). 
\eeq
Here, $X(y)\in\Gamma$ denotes a point near $X(0)$, where $y \in \Rbb^2$ denotes its local coordinate in the chosen chart.
Since $ \det g(y)=1+O(|y|^2)$ in geodesic normal coordinates centered at $x_0=X(0)$, 
the kernel of the single-layer potential operator has the expansion
\beq
\frac{-\sqrt{\det g(y)}}{4\pi |X(y)-X(0)|} = -\frac{1}{4\pi |y|}+O(|y|).
\eeq
It follows that there is no homogeneous term of degree $0$ in the kernel expansion of the single-layer potential operator. Since 
$\Gamma$ is two-dimensional, this implies that the homogeneous component of the symbol of order $-2$
vanishes at $x_0$ in geodesic normal coordinates, that is, $s_{-2}(x_0,\xi)=0.$ 
Here we denote the homogeneous component of the symbol of order $k$ as $s_k(x,\xi)$. Then it follows from the calculus of $\Gy$DOs \cite{Agranovich} that the subprincipal symbol 
$\Scal_\Gamma$ is given by 
\beq
\sub(\Scal_\Gamma)(x,\xi) =
s_{-2}(x,\xi) + \frac{i}{2} \sum_j \partial_{x_j}\partial_{\xi_j}s_{-1}(x,\xi).
\eeq
Moreover, since $\left. \partial_{x_j}g^{kl}(x) \right|_{x=x_0} =0,\ ( j,k,l=1,2)$ in geodesic normal 
coordinates centered at $x_0$, we obtain
\begin{align}
\left.\sub(\Scal_\Gamma)(x,\xi)\right|_{x=x_0}
&= s_{-2}(x_0,\xi) + \frac{i}{2} \sum_j \left. \partial_{x_j}\partial_{\xi_j}
\left[ -\frac{1}{2} \bigl(g^{kl}(x)\xi_k\xi_l\bigr)^{-1/2} \right] \right|_{x=x_0}  \notag\\
&= 0+ \frac{i}{2}\sum_j \partial_{\xi_j}\left[ \frac{1}{4} \bigl(g^{kl}(x)\xi_k\xi_l\bigr)^{-3/2} \left.
\partial_{x_j}g^{kl}(x) \right|_{x=x_0} \xi_k\xi_l \right] \notag\\
&=0. \label{eq: def subprincipal symbol}
\end{align}
Since $x_0\in\Gamma$ is arbitrary, it follows that 
\beq
\sub(\Scal_\Gamma)=0.
\eeq

It is also known that the operator $-\Scal_{\Gamma}$  is  positive and invertible, it is  elliptic, and therefore the complex powers of $-\Scal_\Gamma$ are $\Gy$DOs \cite{Agranovich}. The operator $(-\Scal_{\Gamma})^z$ has principal symbol $(2 |\xi|_g)^{-z}$. Therefore the principal symbol of $\widetilde{\Kcal_{\Gamma}}$ equals the product of principal symbols of $\Kcal_\GG, (-\Scal_\GG)^{1/2}$ and $(-\Scal_\GG)^{-1/2}$, and thus, the principal symbol of the symmetrized NP operator coincides with the one of the original NP operator. In this way, we obtained:
\begin{prop}\label{prop: symbol of symmtrized op}
Let $\Omega \subset{\mathbb R}^3$ be a bounded open set with smooth boundary. Then the single-layer potential operator $\Scal_{\Gamma}$ and the symmetrized NP operator $\widetilde{\Kcal_{\Gamma}}$ are  $\Gy$DOs of order $-1$on $\Gamma$ with principal symbols $p_S =(-2|\xi|_x)^{-1}$ and $p_{\Kcal}(x,\xi)$, respectively. Here $p_{\Kcal}(x,\xi)$ is given by \eqref{eq: NP-PrincSymbol}. Moreover their subprincipal symbols satisfy $\mbox{sub}(\Scal_{\Gamma})=\mbox{sub}(\widetilde{\Kcal_{\Gamma}})=0$.
\end{prop}
Since the vanishing of $\mbox{sub}(\widetilde{\mathcal{K}_\Gamma})$ has not been justified above, we provide a proof of this fact here.
Let us compute the subprincipal symbol of the symmetrized NP operator $\widetilde{\Kcal_{\Gamma}}$, and see that they are at least pure imaginary on local charts. To do so, we recall the calculus of subprincipal symbols.  It is emphasized that the subprincipal symbol is uniquely determined on the half density forms $L^2{(\Omega^{1/2})}$ (H\"olmander's notation \cite{Ho4}).

To simplify the notation, we define $S_0$ and $K_0$ as follows.  
$$
S_0 := \frac{\sqrt{EG-F^2}}{2(G\Gx_1^2 - 2F \Gx_1 \Gx_2 + E\Gx_2^2)^{1/2}}
$$
and
$$
K_0 := -\frac{{\sqrt{EG-F^2}(L \Gx_2^2 - 2M \Gx_1 \Gx_2 + N \Gx_1^2)}}{4 (G\Gx_1^2 - 2F \Gx_1 \Gx_2 + E\Gx_2^2)^{3/2}},
$$
where $E$, $F$, and $G$ are the coefficients of the first fundamental form, and $L$, $M$, and $N$ are those of the second fundamental form, evaluated at $X(x)\in\Gamma$. Thus the notations $S_0$ and $K_0$, respectively, denote the principal symbols for $-\Scal_{\Gamma}$ and $\Kcal_{\Gamma}$. We employ two lemmas to calculate subprincipal symbols \cite{Agranovich}: 
\begin{lem}[\cite{Agranovich}] For $\Gy$DOs $A$ and $B$, we have 
\beq
\sub (AB) = A_0 \sub B + B_0 \sub A + (2i)^{-1} \{ A_0, B_0 \},
\eeq
where $A_0$ and $B_0$ are the principal symbols of $A$ and $B$, respectively, and the Poisson bracket $\{ A_0, B_0\}$ is defined by
\beq
\{ A_0, B_0 \} := \sum_{j = 1}^2 (\p_{\Gx_j} A_0 \p_{x_j} B_0 - \p_{\Gx_j} B_0 \p_{x_j} A_0).
\eeq
\end{lem}

\begin{lem}[\cite{Agranovich}]
For $z \in \Cbb$ and a $\Gy$DO $A$, we have
$$
\sub A^z = z A_0^{z - 1} \sub A.
$$
\end{lem}

Recall that $(A^z)_0=A_0^z$ for all $z\in\Cbb$, where the subscript $0$ denotes the principal symbol
of the corresponding operator. In particular, for any $n\in\Nbb$, we have
\begin{align*}
& \sub(\sqrt{{-\Scal}_{\Gamma}}^{-1} \Kcal_{\Gamma}^{n} \sqrt{-\Scal_{\Gamma} }) \\
=& (\sqrt{-\Scal_{\Gamma}}^{-1})_0 (\Kcal_{\Gamma}^{n})_0 \sub (\sqrt{{-\Scal}_{\Gamma}}) + (\Kcal_{\Gamma}^{n})_0 (\sqrt{{-\Scal}_{\Gamma}})_0 \sub (\sqrt{{-\Scal}_{\Gamma}}^{-1})\\
& + (\sqrt{{-\Scal}_{\Gamma}})_0 (\sqrt{{-\Scal}_{\Gamma}}^{-1})_0 \sub (\Kcal_{\Gamma}^{n}) + (2i)^{-1} (\sqrt{{-\Scal}_{\Gamma}})_0 \{(\sqrt{{-\Scal}_{\Gamma}}^{-1})_0, (\Kcal_{\Gamma}^{n})_0 \}\\
&+ (2i)^{-1} (\sqrt{{-\Scal}_{\Gamma}}^{-1})_0 \left\{ (\Kcal_{\Gamma}^{n})_0, (\sqrt{{-\Scal}_{\Gamma}})_0 \right\}\\
& + (2i)^{-1} (\Kcal_{\Gamma}^{n})_0 \left\{ (\sqrt {{-\Scal}_{\Gamma}}^{-1})_0, (\sqrt{-\Scal_{\Gamma}})_0 \right\}\\
=& S_0^{-1/2} K_0^{n} \frac{1}{2} S_0^{-1/2} \sub (-\Scal_{\Gamma}) + K_0^{n} S_0^{1/2} \left( -\frac{1}{2} S_0^{-3/2} \right) \sub(-\Scal_{\Gamma})\\
&+ S_0^{1/2} S_0^{-1/2} \left( n K_0^{(n - 1)} \right) \sub \Kcal_{\Gamma} + (2i)^{-1} S_0^{1/2} \{ S_0^{-1/2}, K_0^{n} \}\\
&+ (2i)^{-1} S_0^{-1/2} \left\{ K_0^{n}, S_0^{1/2} \right\} + (2i)^{-1} K_0^{n} \left\{ S_0^{-1/2}, S_0^{1/2} \right\}.
\end{align*}
One can easily see that
\begin{align*}
\left\{ S_0^{-1/2}, K_0^{n} \right\} =& -\frac{n}{2} S_0^{-3/2} K_0^{(n - 1)} \{ S_0, K_0 \},\\
\left\{ K_0^{n}, S_0^{1/2} \right\}  =& \frac{n}{2} S_0^{-1/2} K_0^{(n - 1)} \{ K_0, S_0 \} = -\frac{n}{2} S_0^{-1/2} K_0^{(n - 1)} \{ S_0, K_0 \},\\
\left\{ S_0^{-1/2}, S_0^{1/2} \right\} =& -\frac{1}{4} S_0^{-2} \left\{ S_0, S_0 \right\} = 0.
\end{align*}
Thus, we have
\beq\label{subprincipal of Ple NP -n order}
\sub(\sqrt{-\Scal_{\Gamma}}^{-1} \Kcal_{\Gamma}^{n} \sqrt{-\Scal_{\Gamma}}) = n K_0^{(n - 1)} \left( \sub \Kcal_{\Gamma} - (2i)^{-1} S_0^{-1} \{ S_0, K_0 \} \right).
\eeq

We compute the sub-principal symbol of $\Kcal_{\Gamma}$. It is given by
$$
\sub \Kcal_{\Gamma} = k_{-2} + \frac{i}{2} \sum_{|\Ga| = 1} \p_{\Gx}^{\Ga} \p_x^{\Ga} K_0,
$$
where $k_{-2}$ is the symbol of $\Kcal$ of order $-2$. We find that $k_{-2}(x, \xi)$ is purely imaginary. In fact, the local expansion of the integral kernel of $\Kcal_{\Gamma}$ shows that the homogeneous term of degree zero is real-valued and odd in the difference variable. 
Its Fourier transform is therefore purely imaginary, which implies that $k_{-2}(x, \xi)$ is purely imaginary. Thus $k_{-2}$ is pure imaginary on local charts although we do not have its explicit formula. Since $K_0$ is real, the sub-principal symbol $\sub \Kcal_{\Gamma}$ is purely imaginary.
We can also find that $(2i)^{-1} S_0^{-1} \{ S_0, K_0 \}$ in \eqref{subprincipal of Ple NP -n order} is pure imaginary, which implies that the sub-principal symbol 
$\sub(\sqrt{-\Scal_{\Gamma}}^{-1} \Kcal_{\Gamma}^{n} \sqrt{-\Scal_{\Gamma}})$ is also {\it pure imaginary}.  

Furthermore $\sqrt{-\Scal_{\Gamma}}^{-1} \Kcal_{\Gamma}^{n} \sqrt{-\Scal_{\Gamma}}$ is self-adjoint in $L^2(\Gamma)$ by Plemelj's principle. 

On the other hand, it is known that, for an $m$-th order polyhomogeneous pseudodifferential operator $A$ satisfying
$A-A^{*}\in \Psi_{\rm ph}^{m-2},$ we have $A_0=\overline{A_0}, \ \mbox{sub}A=\overline{\mbox{sub}A}.$
In particular, the principal and subprincipal symbols of the self-adjoint operator $A$ are real-valued. It therefore follows that
\beq
\mbox{sub}( \sqrt{-\Scal_{\Gamma}}^{-1} \Kcal_{\Gamma}^{n} \sqrt{-\Scal_{\Gamma}} )
\eeq
is {\it real-valued}. Thus we conclude that
\beq
\sub(\widetilde{\Kcal_{\Gamma}}^{n})= 
\sub(\sqrt{-\Scal_{\Gamma}}^{-1} \Kcal_{\Gamma}^{n} \sqrt{-\Scal_{\Gamma}}) = 0.
\eeq
Especially when $n=1$, $\sub(\widetilde{\Kcal_{\Gamma}})= 0$.

\subsection{The Neumann--Poincar\'e flow}\label{sec: NP flow}
As is well known \cite{Ho4}, in order to obtain precise eigenvalue asymptotics of the first order $\Gy$DOs, one first needs to consider the Hamiltonian flow generated by the principal symbol of order $1$. However, since the single-layer potential operator and the NP operator are pseudodifferential operators of order $-1$, it is necessary to consider their inverse $\Scal_{\Gamma}^{-1}$ and $\Kcal_{\Gamma}^{-1}$, which are $\Gy$DOs of order $1$, in order to carry out a more detailed analysis. If we consider the principal symbol of $\Scal_{\Gamma}^{-1}$ as a Hamiltonian on $T^{*}\Gamma$, then the corresponding Hamiltonian flow is just the geodesic flow on $S^{*}\Gamma$. Here $T^{*}\Gamma$ and $S^{*}\Gamma$ respectively denote the cotangent bundle and the cosphere bundle.  
Similarly, for the case of $\Kcal_{\Gamma}^{-1}$, we refer to the corresponding flow as the Neumann--Poincar\'e (NP) flow.
\begin{definition}[The Neumann--Poincar\'e flow]\label{def: NP Hamilton flow}
Let $\GO \subset \Rbb^3$ be a bounded open set with smooth boundary $\Gamma=\p\GO$, and assume that
the Gaussian curvature of $\Gamma$ is everywhere positive. We define $p(x, \xi):= \frac{1}{|K_0(x, \xi)|}$ for $\xi \not =0$. Then the NP energy surface $\Sigma$, the NP Hamiltonian vector field and the corresponding NP Hamilton flow are, respectively, defined by  
$$
\begin{cases}
&\Sigma=\{ (x, \xi)\in T^{*}\Gamma \setminus \{\xi =0\}\ |\ p(x, \xi)=1 \}, \\
&X_{p}:=\displaystyle\sum_{i=1}^{2} \left(\frac{\partial p}{\partial \xi_i}\frac{\partial}{\partial x_i}- 
\frac{\partial p}{\partial x_i}\frac{\partial}{\partial \xi_i} \right), \\
&\exp \tau X_{p}: \Sigma\ \rightarrow \ \Sigma.
\end{cases}
$$
We define the Liouville measure $d\nu$ on $\Sigma$ by the Gelfand--Leray formula
\beq
d\nu = \left. \frac{dx_1\wedge d\xi_1 \wedge dx_2\wedge d\xi_2}{dp} \right|_{\Sigma}.
\eeq
 Here the measure $d\nu$ is invariant under the Hamiltonian flow generated by $p$ (See, e.g., \cite{Ar, D} for the details of the Liouville measure.).  
\end{definition}
We emphasize that the NP Hamiltonian flow is essential to the second term of Weyl's law of the NP operator. 

\begin{definition}[The NP periodic points]\label{def: NP periodic points}
In the same setting as in Definition \ref{def: NP Hamilton flow}, we define the set of the NP periodic points by 
\beq
S_{\mbox{\rm periodic}} := \{ (x, \xi) \in \Sigma\ |\ \exists \tau \neq 0 \ \mbox{such that}\ \exp \tau X_{p} (x, \xi) =(x, \xi)\}.
\eeq
If $\nu(S_{\mbox{\rm periodic}})=0$, we say that the set of periodic points of the NP flow has measure zero.
\end{definition}
For the boundary of finite smoothness $C^{2,\alpha}$, some of the above definitions remain valid. 
However, since a convenient  symbolic calculus for operators on such manifolds is not developed, we restrict our attention to the smooth case.

\subsection{Notes on eigenvalues of  the Dirichlet--to--Neumann map}\label{sec: DN map eigenvalues notation}
Although we restrict our attention to smooth boundaries in this paper, we briefly mention the relevant results for the asymptotic behaviors of eigenvalues of the Dirichlet--to--Neumann (abbreviated by DtN) map $\Lambda$. We begin by recalling the following general result:  
\beq 
\lambda_j(\Lambda)=2 \pi\left(\frac{j}{\text{Vol}(B^{d-1})\text{Vol}(\Gamma)}\right)^{\frac{1}{d-1}}+O(1)
\eeq
where $B^{d-1}$ is the unit ball in $\Rbb^{d-1}$.

We denote by $\Lambda^{-1}$ the inverse of $\Lambda$ on the orthogonal complement of $\ker\Lambda$,
extended by zero on $\ker\Lambda$. Then the nonzero eigenvalues of $\Lambda^{-1}$, which is also known as 
the Neumann--to--Dirichlet (NtD) map on the appropriate subspace, satisfy
\beq\label{eq: Neumann-Poincare eigenvalue asymptotic in higher dimension}
\lambda_j(\Lambda^{-1})=\frac{1}{2\pi} \left(\frac{\text{Vol}(B^{d-1})\text{Vol}(\Gamma)}{j}\right)^{\frac{1}{d-1}}+O(j^{\frac{-2}{d-1}}).
\eeq

We also note that more detailed results are available in the two-dimensional case \cite{GPPS}.  
To state this precisely, let  $C=(\alpha_1, \alpha_2, \ldots, \alpha_k)$ be a finite sequence of positive numbers, where repetitions are allowed. We define the sets of multiplicities by 
\beq
S(C)=
\{\underbrace{0,\ldots,0}_{k \text{\ zeros}}\}
\biguplus_{i=1}^k
\left(
\alpha_i\mathbb N_{\ge1}
\biguplus
\alpha_i\mathbb N_{\ge1}
\right).
\eeq
Here, the notation $\biguplus$ denotes the multiset union, 
where the multiplicity of each element is the sum
of its multiplicities in the individual multisets. For instance, 

$S(\{1\}) = \{0, 1, 1, 2, 2, 3, 3, \cdots \}$ and $S(\{2, \pi \})=\{0, 0, 2, 2, \pi, \pi, 4, 4, 6, 6, 2\pi, 2\pi, \cdots \}$.

Then the precise asymptotic formula is given in the following theorem.  
\begin{thm}[A. Girouard, L. Parnovski, I. Polterovich, D. Sher (2014) \cite{GP, GPPS} ]\label{thm: DN map eigenvalues}
Let $\Gamma_1, \cdots \Gamma_k$ be the connected components of the boundary $\Gamma=\p\GO \subset \Rbb^2$, with lengths 
$l(\Gamma_i),\ 1\leq i \leq k$. We define $C$ by the finite sequence of positive numbers: 
\beq
C= \left( \frac{2\pi}{l(\Gamma_1)}, \cdots, \frac{2\pi}{l(\Gamma_k)} \right).
\eeq
Then
$\lambda_j(\Lambda) = S(C)_j + O(j^{-\infty}). $ 
Here the notation $S(C)_j$ denotes the $j$-th smallest element of $S(C)$, counted with multiplicity.The notation $O(j^{-\infty})$ means the error term decays faster than any negative power of $j$.
\end{thm}
This theorem also plays an important role in the spectral analysis of the single-layer potential operator in Section \ref{sec: single-layer spectrum 2D}. We emphasize that the $k$ zeros in $S(C)$ should not be confused with either the dimension of the kernel of the DtN map or that of the single-layer potential operator $\mathcal S_\Gamma$. The number $k$ denotes only the number of connected components of $\Gamma$. 

For the DtN map, Theorem \ref{thm: DN map eigenvalues} is formulated in a general setting that includes the case of a single boundary component as a special case: 
\begin{cor}[\cite{E, R}]
Assume that the region $\GO$  has a connected boundary $\Gamma=\p\Omega$, with
lengths $l(\Gamma)$. Then 
\beq
\lambda_{2j} (\Lambda)=\lambda_{2j+1}(\Lambda) +O(j^{-\infty}) =\frac{2\pi}{l(\Gamma)}j +O(j^{-\infty}) \ \text{as}\ j \to \infty.
\eeq
Here $O(j^{-\infty})$ means that the error term decays faster than any negative power of j.
\end{cor}

\section{Eigenvalue asymptotics of layer potential operators for $C^{\infty}$ boundaries}\label{sec: Eigenvalue asymptotics of layer potential operators for C^infty}
In this section, we investigate the asymptotic behavior of the eigenvalues of the single-layer potential operator and the Neumann–Poincaré (NP) operator.  We first study the single-layer potential operator in two and three dimensions in separate subsections.
\subsection{Eigenvalues of the single-layer potential operator in two dimensions}\label{sec: single-layer spectrum 2D}

\textit{(Outline of the proof and interpretation of $\Scal_\Gamma^{-1}$).}
Since $\Scal_\Gamma$ may have a nontrivial kernel, we introduce its Moore--Penrose pseudoinverse
$\Scal_\Gamma^\dagger$. Let $\Pi_0$ be the orthogonal projection onto $\ker\Scal_\Gamma$, and we define
$\Scal_\Gamma^\dagger$ as the inverse of $\Scal_\Gamma$ on $(\ker\Scal_\Gamma)^\perp$, extended  by zero on the kernel. 
Then the operator $B=-\Scal_\Gamma^\dagger$ is a self-adjoint elliptic pseudodifferential operator of order one 
with positive principal symbol $2|\xi|_g$. As shown below, the jump relation implies
\beq
B-2\Lambda\in\Psi^{-\infty}(\Gamma).
\eeq
The spectral stability under smoothing
perturbations, together with
Theorem \ref{thm: DN map eigenvalues},
then yields the desired eigenvalue asymptotics.
The possible kernel of $\Scal_\Gamma$ is
accounted for by a shift in the eigenvalue indices.

\begin{pf}[Proof of Theorem
\ref{thm: single-layer accurate 2D eigenvalues}]

Recall that the jump relation
\eqref{eq: relation among single, DN and NP}
gives
\beq\label{eq: Jump identity}
\Scal_\Gamma\Lambda
=\Kcal_\Gamma-\frac{1}{2} I.
\eeq
We multiply both sides of \eqref{eq: Jump identity} by $\Scal_\Gamma^\dagger$ from the left. 
It then follows from $\Scal_\Gamma^\dagger\Scal_\Gamma=I-\Pi_0$ that 
\beq
(I-\Pi_0)\Lambda
=\Scal_\Gamma^\dagger\Kcal_\Gamma
-\frac12\Scal_\Gamma^\dagger.
\eeq
Consequently,
\beq\label{eq: pseudoinverse DN relation}
\Scal_\Gamma^\dagger+2\Lambda
=2\Scal_\Gamma^\dagger\Kcal_\Gamma
+2\Pi_0\Lambda.
\eeq

It is known that $\Scal_\Gamma$ is an elliptic $\Gy$DO of order $-1$ (e.g., \cite{PS}).
Hence its pseudoinverse $\Scal_\Gamma^\dagger$ is an elliptic $\Gy$DO of order one.

Moreover, for smooth planar boundaries, $\Kcal_\Gamma$ is a smoothing operator (e.g., \cite{MS}).
It follows from the elliptic regularity of $\Scal_\Gamma$ that every element of $\ker\Scal_\Gamma$ is smooth.
Since $\dim\ker \Scal_{\Gamma}<+\infty$, $\Pi_0$ and $\Pi_0\Lambda$ are smoothing operators.

Thus it follows from \eqref{eq: pseudoinverse DN relation} that
\beq
B-2\Lambda =-2\Scal_\Gamma^\dagger\Kcal_\Gamma -2\Pi_0\Lambda \in\Psi^{-\infty}(\Gamma).
\eeq

We now apply the standard spectral stability result for self-adjoint elliptic
$\Gy$DOs of positive order under smoothing perturbations (see, e.g., \cite[Lemma 2.1]{GPPS}).

Let $\mu_j(B)$ denote the eigenvalues of $B$
arranged in nondecreasing order and counted
with multiplicity.
Since $B$ is bounded from below, we may add
a sufficiently large constant $c>0$ so that
both $B+cI$ and $2\Lambda+cI$ are positive.
Applying the spectral stability result to
these shifted operators gives
\beq
\mu_j(B)
=\lambda_j(2\Lambda)+O(j^{-\infty}).
\eeq
Here we know from Theorem \ref{thm: DN map eigenvalues} that 
\beq
\lambda_j(2\Lambda) =4\pi S(R)_j+O(j^{-\infty}),
\eeq
where $R$ is the sequence of reciprocals of the lengths of the connected components
of $\Gamma$, as defined in Theorem \ref{thm: single-layer accurate 2D eigenvalues}. Thus,
\beq
\mu_j(B)
=4\pi S(R)_j+O(j^{-\infty}).
\eeq

Finally, let
$m_0=\dim\ker\Scal_\Gamma$.
Since $B=-\Scal_\Gamma^\dagger$,
the nonzero eigenvalues of $\Scal_\Gamma$
are the negative reciprocals of the
nonzero eigenvalues of $B$.
Taking into account the $m_0$ zero eigenvalues
of $B$, we obtain, for sufficiently large $j$,
\beq
\lambda_j(\Scal_\Gamma)
=-\frac{1}{\mu_{j+m_0}(B)}.
\eeq
Since $S(R)_j$ grows linearly in $j$,
taking reciprocals preserves the
$O(j^{-\infty})$ remainder.
Therefore,
\beq
\lambda_j(\Scal_\Gamma)
=-\frac{1}{4\pi}
\left(S(R)_{j+\dim\ker\Scal_\Gamma}\right)^{-1}
+O(j^{-\infty}).
\eeq
This completes the proof.
\end{pf}
\subsubsection{Further consequences of Theorem 1.1}
Even for smooth boundaries $\Gamma$ of planar domains, one might expect a second-order term of the form
$cj^{-2}$, where $c$ is a constant, following the leading term of order $j^{-1}$.
However, although the remainder is of order $O(j^{-2})$, such a constant second-order coefficient does not exist in general. 

As we show below, the asymptotic pairing of successive eigenvalues gives rise to oscillations of order $j^{-2}$. 
This phenomenon is consistent with the fact that all trajectories of the Hamiltonian flow generated by $p(x,\xi)=2|\xi|_g$ 
are periodic on each connected component of the one-dimensional boundary $\Gamma$.
Consequently, the nonperiodicity assumption employed in Section \ref{sec:  single-layer potentials in three dimensions}
to improve the remainder estimate in three dimensions is never satisfied in the planar case.

More generally, when the set of periodic points of the Hamiltonian flow has positive Liouville measure,
the improved Weyl remainder estimate is no longer guaranteed. In such cases, periodic trajectories
may give rise to oscillatory contributions to the Weyl remainder or to the formation of eigenvalue clusters.
We refer to \cite{DG,I,Sa} for related results on spectral asymptotics.

In fact, the asymptotic pairing of successive eigenvalues yields the oscillations of order
$j^{-2}$ in two dimensions, such a behavior prevents the existence of the second-order asymptotic term 
with a constant coefficient. Our porpose in this subsection is to show this eigenvalue behavior precisely.

To show this fact, let $\Gamma_1,\ldots,\Gamma_k$ be the connected
components of $\Gamma$, and let $l_i=l(\Gamma_i)$. We first remark that 
the number of zero eigenvalues of the DtN map need not coincide with 
the number of zero entries in the model sequence $S(C)$. In fact, if 
$\Omega$ has $q$ connected components, then 
\beq
\dim\ker\Lambda=q,
\eeq
whereas $S(C)$ contains $k$ zeros.
These zeros are part of the indexing convention in Theorem \ref{thm: DN map eigenvalues} and are essential
for the $O(j^{-\infty})$ remainder estimate.

\begin{cor}
Under the assumptions of Theorem \ref{thm: single-layer accurate 2D eigenvalues}, we have the followings: 

Let $\displaystyle L=l(\Gamma)=\sum_{i=1}^k l_i$. Then
\beq
\lambda_j(\Scal_\Gamma)=-\frac{L}{2\pi j}+O(j^{-2}). 
\eeq
Moreover, 
\beq\label{eq: non-existence of the second term of eigenvalue asymptotics for the single layer potential operator}
\limsup_{j\to\infty}
j^2\left(\lambda_j(S_\Gamma)+\frac{L}{2\pi j}\right) - \liminf_{j\to\infty}
j^2\left(\lambda_j(S_\Gamma)+\frac{L}{2\pi j}\right) \geq \frac{L}{2\pi} 
\eeq
In particular, there is no constant $c$
such that
\beq
\lambda_j(S_\Gamma)
=-\frac{L}{2\pi j}
+\frac{c}{j^2}+o(j^{-2}).
\eeq
\end{cor}

\begin{pf}

Let $k$ denote the number of connected components of $\Gamma$, and 
\beq
m_0:=\dim\ker S_\Gamma,\qquad
l_i:=l(\Gamma_i),\qquad
L:=\sum_{i=1}^k l_i.
\eeq

Remarking that the sequence $S(R)$ contains exactly $k$ zeros, counted with multiplicity, we define
\beq
a_j:=\frac{1}{4\pi S(R)_{j+k}},
\quad j\geq 1.
\eeq

It follows from the definition of $S(R)$ that the sequence
$\{a_j\}$ is the nonincreasing rearrangement counting multiplicities of the multiset
\beq
\biguplus_{i=1}^k
\left\{
\frac{l_i}{4\pi n},
\frac{l_i}{4\pi n}:n\geq 1
\right\}.
\eeq
In particular, $a_{2m-1}=a_{2m}, m\geq 1.$ 

It follows from Theorem \ref{thm: single-layer accurate 2D eigenvalues} that we have for all sufficiently large $j$,
\beq
\lambda_{j+k-m_0}(S_\Gamma)
=-a_j+O(j^{-\infty}).
\eeq

The counting function associated with the sequence
$\{a_j\}_{j\geq 1}$ satisfies
\begin{align}
N(t)
&:=\#\{j\geq 1:a_j>t\} \notag\\
&=2\sum_{i=1}^k
\#\left\{n\geq 1:
\frac{l_i}{4\pi n}>t\right\} \notag\\
&=\frac{L}{2\pi t}+O(1),
\quad t\downarrow 0.
\end{align}
Since the multiplicities in the sequence are uniformly bounded, it follows from inverting the counting function that  
\beq\label{eq: a_j estimates}
a_j=\frac{L}{2\pi j}+O(j^{-2}).
\eeq
Since $k-m_0$ is a constant, we obtain
\beq\label{eq: rough single layer potential asymtotic formula}
\lambda_j(S_\Gamma)=-\frac{L}{2\pi j}+O(j^{-2}).
\eeq

We next show that the remainder of order $j^{-2}$ does not admit a constant coefficient.

For each $m$, we define $p_m=2m+k-m_0-1,\ q_m=2m+k-m_0=p_m+1.$
It then follows from Theorem \ref{thm: single-layer accurate 2D eigenvalues} and the identity $a_{2m-1}=a_{2m}$ that 
\beq
\lambda_{p_m}(S_\Gamma)-\lambda_{q_m}(S_\Gamma)=O(m^{-\infty}).
\eeq

It follows from \eqref{eq: rough single layer potential asymtotic formula} that the sequence
\beq
E_j:=j^2\left(\lambda_j(S_\Gamma)+\frac{L}{2\pi j}\right)
\eeq
is bounded. 
Moreover,
\begin{align*}
E_{q_m}-E_{p_m}
&=q_m^2\lambda_{q_m}(S_\Gamma)
-p_m^2\lambda_{p_m}(S_\Gamma)+\frac{L}{2\pi}\\
&=(q_m^2-p_m^2)\lambda_{q_m}(S_\Gamma)
+p_m^2\bigl(
\lambda_{q_m}(S_\Gamma)
-\lambda_{p_m}(S_\Gamma)
\bigr)+\frac{L}{2\pi}.
\end{align*}
Since
$
\lambda_{q_m}(S_\Gamma)
=-\frac{L}{2\pi q_m}+O(q_m^{-2})
$
and $q_m-p_m=1$, we conclude that
\beq
\lim_{m\to\infty}
(E_{q_m}-E_{p_m})=-\frac{L}{2\pi}.
\eeq

Consequently,
\beq
\limsup_{j\to\infty}E_j
-\liminf_{j\to\infty}E_j
\geq \frac{L}{2\pi}>0.
\eeq

This proves the inequality \eqref{eq: non-existence of the second term of eigenvalue asymptotics for the single layer potential operator}.

In particular, there is no constant $c$
such that
$$
\lambda_j(S_\Gamma)
=-\frac{L}{2\pi j}
+\frac{c}{j^2}+o(j^{-2}).
$$

\end{pf}

\subsubsection{Examples of eigenvalues of the single-layer potential operator in two dimensions}\label{sec: Examples of eigenvalues of the single-layer potential operator }
We present two classical examples, a disk and an annulus, for which the eigenvalues can be computed explicitly, thereby illustrating Theorem \ref{thm: single-layer accurate 2D eigenvalues}.  

Let $\GO=B(0,R)$ be a disk of radius $R$, and let $\GG=\p\GO$.
For $n\in\Zbb$, define the density $\Gvf_n(y)=e^{in\Gt_y}$, where $y=R(\cos\Gt_y,\sin\Gt_y)\in\p\GO$.
In polar coordinates $x=r(\cos\Gt,\sin\Gt)$, the following formulas are well known (see, e.g., \cite{Z}):
\begin{align}
\label{slayer_disk}
\Scal_{\GO}[\Gvf_n](x) &= \begin{cases}
\displaystyle -\frac{R}{2|n|}\left(\frac{r}{R}\right)^{|n|}e^{in\Gt}, & r\leq R, \\
\displaystyle -\frac{R}{2|n|}\left(\frac{R}{r}\right)^{|n|}e^{in\Gt}, & r>R,
\end{cases}
\quad n\in\Zbb_{\neq 0},
\end{align}
and
\begin{align}
\label{slayer_disk_constant}
\Scal_{\GO}[\Gvf_0](x) &= \begin{cases}
R\ln R, & r\leq R, \\
R\ln r, & r>R.
\end{cases}
\end{align}
Since $\{\Gvf_n\}_{n\in\Zbb}$ is a complete orthogonal basis of $L^2(\GG)$, the eigenvalues of $\Scal_\GG$ are 
\begin{align}
R \log R \text{ (multiplicity 1)}, \quad -\frac{R}{2n} \text{ (multiplicity 2)}, \; n \geq 1.
\end{align}

If $R\neq1$, then $\ker\Scal_\GG=\{0\}$, and the nonzero
eigenvalue $R\log R$ is included in the ordered sequence.
Thus, for all sufficiently large $j$,
\begin{align}
\Gl_{2j}(\Scal_\GG)=\Gl_{2j+1}(\Scal_\GG)
=-\frac{R}{2j}=-\frac{l(\GG)}{4\pi j}.
\end{align}
If $R=1$, then $\dim\ker\Scal_\GG=1$, and
\begin{align}
\Gl_{2j-1}(\Scal_\GG)=\Gl_{2j}(\Scal_\GG)
=-\frac{R}{2j}=-\frac{l(\GG)}{4\pi j},
\qquad j\geq1.
\end{align}
These identities verify
\eqref{eq: single layer dimker Scal=0 case} and
\eqref{eq: single layer dimker Scal=1 case}, respectively,
in Theorem \ref{thm: single-layer accurate 2D eigenvalues}.

We next consider the annulus $\GO = B(0,R_2) \setminus \overline{B(0,R_1)}$, where $0<R_1<R_2$, and set $\Gr:=R_1/R_2$ ($0<\Gr<1$). Its boundary $\GG=\p\GO$ consists of two connected components, the inner circle $\GG_1$ and the outer circle $\GG_2$. We denote the angular densities on each boundary circle by $f_n(y)=e^{in\Gt_y}$, $n\in\Zbb$. We identify a density $\Gvf$ on $\GG$ with the column vector
\begin{align}
\label{eq:boundary_density_identification}
\Gvf &\longleftrightarrow \begin{bmatrix}\Gvf_1 \\ \Gvf_2\end{bmatrix},
\end{align}
where $\Gvf_1:=\Gvf\big|_{\GG_1}$ and $\Gvf_2:=\Gvf\big|_{\GG_2}$.
It then follows from \eqref{slayer_disk} and \eqref{slayer_disk_constant} that
\begin{align}
\Scal_{\GG} \left[\begin{array}{c}{f_n} \\ {0}  \end{array}\right](x) &= -\frac{1}{2|n|}\begin{bmatrix}R_1 e^{in\Gt} \\ R_1\Gr^{|n|}e^{in\Gt}\end{bmatrix}, \quad n \in \Zbb_{\neq 0}, \\ 
\Scal_{\GG} \left[\begin{array}{c}{0} \\ {f_n}  \end{array}\right](x) &= -\frac{1}{2|n|}\begin{bmatrix}R_2\Gr^{|n|}e^{in\Gt} \\ R_2 e^{in\Gt}\end{bmatrix}, \quad n \in \Zbb_{\neq 0}, 
\end{align}
and
\begin{align}
\Scal_{\GG} \left[\begin{array}{c}{f_0} \\ {0}  \end{array}\right](x) &= \left[\begin{array}{c}{R_1 \ln R_1} \\ {R_1 \ln R_2}  \end{array}\right], \\
\Scal_{\GG} \left[\begin{array}{c}{0} \\ {f_0}  \end{array}\right](x) &= \left[\begin{array}{c}{R_2 \ln R_2} \\ {R_2 \ln R_2 }  \end{array}\right].
\end{align}
Therefore, for each $n\in\Zbb$, the subspace spanned by $[f_n,0]^T$ and $[0,f_n]^T$ is invariant under $\Scal_{\GG}$.
With respect to this ordered basis, the restrictions of $\Scal_{\GG}$ are represented by
\begin{align}
\label{eq:annulus_nonzero_matrix}
A_n &= -\frac{1}{2|n|}
\begin{bmatrix}
R_1 & R_2\Gr^{|n|} \\
R_1\Gr^{|n|} & R_2
\end{bmatrix}, \quad n\in\Zbb_{\neq 0}, \\
\label{eq:annulus_constant_matrix}
A_0 &= \begin{bmatrix}
R_1\ln R_1 & R_2\ln R_2 \\
R_1\ln R_2 & R_2\ln R_2
\end{bmatrix},
\end{align}
where the corresponding eigenvalues are given by
\begin{gather}
-\frac{R_2+R_1 \pm \sqrt{(R_2-R_1)^2+4R_1R_2\Gr^{2|n|}}}{4|n|}, \\
\frac{R_1\ln R_1+R_2\ln R_2}{2} \pm \frac{1}{2}\sqrt{(R_2\ln R_2-R_1\ln R_1)^2+4R_1R_2(\ln R_2)^2}.
\end{gather}

Since $A_n=A_{-n}$, we define $\Gm_{n}^{(i)}$ by
\begin{align}
\Gm_{2n-1}^{(i)}=\Gm_{2n}^{(i)}
&=-\frac{R_2+R_1+(-1)^i\sqrt{(R_2-R_1)^2+4R_1R_2\Gr^{2n}}}{4n},\quad n\in\Nbb_{\ge1}, \\
\Gm_0^{(i)}
&=\frac{R_1\ln R_1+R_2\ln R_2}{2}+\frac{(-1)^i}{2}\sqrt{(R_2\ln R_2-R_1\ln R_1)^2+4R_1R_2(\ln R_2)^2}.
\end{align}
Since $\{[f_n,0]^T,[0,f_n]^T:n\in\Zbb\}$ is a complete orthogonal basis of $L^2(\GG)$, the numbers $\Gm_n^{(i)}$,
$i=1,2$, $n\geq0$, give all eigenvalues of $\Scal_\GG$,
counting multiplicities. Moreover,
\begin{align}
\Gm_{2n-1}^{(1)}=\Gm_{2n}^{(1)}
&=-\frac{R_1}{2n}+O\left(\frac{\Gr^{2n}}{n}\right)
=-\frac{l(\GG_1)}{4\pi n}+O\left(\frac{\Gr^{2n}}{n}\right),
\label{eq: annulus first family asymptotics}\\
\Gm_{2n-1}^{(2)}=\Gm_{2n}^{(2)}
&=-\frac{R_2}{2n}+O\left(\frac{\Gr^{2n}}{n}\right)
=-\frac{l(\GG_2)}{4\pi n}+O\left(\frac{\Gr^{2n}}{n}\right).
\label{eq: annulus second family asymptotics}
\end{align}
The eigenvalue sequence $\{\Gl_j(\Scal_\GG)\}_{j\geq1}$
is obtained by taking the multiset union of
$\{\Gm_n^{(1)}\}_{n\geq0}$ and $\{\Gm_n^{(2)}\}_{n\geq0}$,
omitting zero entries, and arranging the remaining elements
in nonincreasing order of their absolute values,
counting multiplicities. This ordering agrees with the asymptotic formula in Theorem \ref{thm: single-layer accurate 2D eigenvalues}.

\subsection{Eigenvalues of the single-layer potential operator in three dimensions}\label{sec:  single-layer potentials in three dimensions}
In three and higher dimensions, to the best of our knowledge, the spectral behavior is not as well understood as in the two-dimensional case. Therefore, we establish Theorem \ref{thm: Weyl asymptotic of single-layer potentials in three dimensions}  with a remainder of order $O(j^{-1})$ or $o(j^{-1})$ under the respective assumptions.

\begin{pf}[Proof of Theorem \ref{thm: Weyl asymptotic of single-layer potentials in three dimensions}]
It follows from \eqref{eq: relation among single, DN and NP}  and \eqref{eq: NP-PrincSymbol} that 
\beq 
\Scal_{\Gamma} \equiv -\frac{1}{2} \Lambda^{-1} \quad \text{modulo $\Gy$DO of order }-2.
\eeq
Here $\Lambda^{-1}$ is understood as the inverse on the orthogonal complement of $\ker\Lambda$.
Thus the first asymptotic behavior of eigenvalues of the single-layer potential coincides with one of $-\frac{1}{2} \Lambda^{-1}$, 
that is, 
$$
\lambda_{j}(\Scal_{\Gamma})=\frac{-1}{4\pi} \left(\frac{\text{Vol}(B^{2})\text{Area}(\Gamma)}{j}\right)^{\frac{1}{2}}+O(j^{-1}) = -\frac{1}{4\sqrt{\pi}} {\text{\rm Area}(\Gamma)}^{1/2}  j^{-1/2}+O(j^{-1}) 
$$
as desired. 
Here we can give an alternative explanation for the formula \eqref{eq: second term of single-layer}.  
As mentioned in Section \ref{sec: symbol}, the principal symbol of the single-layer potential operator is $(-2 |\xi|_{x'})^{-1}$.
Then the formula for the singular values $s_j(\Scal_{\Gamma})$ of the single-layer potentials have the similar asymptotics (See, e.g., \cite{BS, Gr, Miyanishi:Weyl} ), that is, 
\beq\label{eq: singular value asymptotic for single-layer}
s_{j}(\Scal_{\Gamma}) = \frac{1}{4\sqrt{\pi}} {\text{\rm Area}(\Gamma)}^{1/2}  j^{-1/2}+O(j^{-1}) .
\eeq
Since $\Scal_\Gamma$ is a negative self-adjoint operator, its eigenvalues are the negatives of its singular values. Consequently,   
\eqref{eq: singular value asymptotic for single-layer} yields \eqref{eq: single-layer rough spectral asymptotic of eigenvalues in 3 dimension}.

To prove the formula \eqref{eq: second term of single-layer}, we use $\sub \Scal_{\Gamma} = 0$, where $\sub \Scal_{\Gamma}$ is the sub-principal symbol of the single-layer potential operator $\Scal_{\Gamma}$. Then the inverse $\Scal_{\Gamma}^{-1}$ is a $\Gy$DO of order $1$, $\sigma_{\rm prin}(\Scal^{-1}_{\Gamma})=-2|\xi|_x$ and $\sub(\Scal_{\Gamma}^{-1})=- \sigma_{\rm prin}(\Scal^{-1}_{\Gamma})^2 \sub(\Scal_{\Gamma})=0$. 
Then the Hamiltonian flow associated with the principal symbol coincides with the geodesic flow.

Here, let us recall a deep result concerning Weyl's law for the first order $\Gy$DOs.
\begin{prop}[The second term of Weyl's law for $\Gy$DOs of order $1$]\label{thm: Weyl law for psiDO of order 1}
Let $M$ be an $n$-dimensional compact smooth manifold without boundary and $A$ be a positive self-adjoint elliptic $\Gy$DO of order $1$ on $M$ with principal symbol $\sigma_{\rm prin}(A)$. Assume that $\mbox{sub}(A)=0$ and the set of periodic points of the Hamiltonian flow generated by $\sigma_{\rm prin}(A)$ has measure zero in $\{ (x, \xi) \in T^{*} M\ |\ \sigma_{\rm prin}(A)(x, \xi)=1\}$. Then 
\beq 
N(A, \lambda)=\alpha \lambda^{n}+o(\lambda^{n-1})\quad \mbox{as} \ \lambda \rightarrow \infty.
\eeq
Here  $N(A, \lambda):=\# \{j | \lambda_j(A) \leq \lambda \}$ is the counting function of the $\Gy$DO $A$ and  $\alpha = (2\pi)^{-n} \mbox{\rm Vol}\{(x, \xi) \in T^{*}M |\ \sigma_{\rm prin}(A)(x, \xi)<1\}. $
\end{prop}

This theorem  for Laplacian was originally conjectured by H. Weyl \cite{Weyl1913}. The first step toward proving Weyl's conjecture for positive elliptic pseudodifferential operators of positive order was made by J. J. Duistermaat and V. W. Guillemin \cite{DG}, under the assumption that the set of periodic points of the associated Hamiltonian flow has Liouville measure zero (see also \cite{Ho4}, \S 29). We also note that estimates for the case in which the set of periodic trajectories has positive measure can be found in H\"ormander's book.

It follows from $\sigma_{\rm prin}(-\Scal_{\Gamma}^{-1})(x, \xi) = 2 |\xi|_x$ that
\begin{align}
\alpha &= (2\pi)^{-n} \mbox{Vol}\{(x, \xi) \in T^{*}\Gamma\ |\ \sigma_{\rm prin}(-\Scal_{\Gamma}^{-1})(x, \xi)<1\} \notag\\
&= (2\pi)^{-2} \mbox{Vol}\{(x, \xi) \in T^{*}\Gamma\ |\ 2 |\xi|_x <1\} \notag\\
&= (2\pi)^{-2} \frac{\pi}{4}\mbox{Area}(\Gamma)  =\frac{\mbox{Area}(\Gamma)}{16 \pi}.
\end{align}
Thus 
\beq
N(-\Scal_{\Gamma}^{-1}, \lambda) = \frac{\mbox{Area}(\Gamma)}{16 \pi} \lambda^2 + o(\lambda) \quad \mbox{as} \ \lambda \rightarrow \infty.
\eeq
For the $j$-th eigenvalue of $-\Scal_{\Gamma}^{-1}$, we have
\beq 
j = \alpha \lambda_j(-\Scal_{\Gamma}^{-1})^2 + o (\lambda_j(-\Scal_{\Gamma}^{-1})),
\eeq
where $\alpha =\frac{{\rm Area}(\Gamma)}{16\pi}.$
Hence, 
\beq 
\lambda_j(-\Scal_{\Gamma}^{-1}) =
\sqrt{\frac{j}{\alpha}} + o(1).
\eeq
It follows that the $j$-th eigenvalue of $-\Scal_{\Gamma}$ satisfies
\beq 
\lambda_j(-\Scal_{\Gamma}) = \frac{1}{\lambda_j(-\Scal_{\Gamma}^{-1})}
= \sqrt{\frac{\alpha}{j}} + o(j^{-1}).
\eeq
Thus,
\beq
\lambda_j(-\Scal_{\Gamma}) =\sqrt{\frac{\mbox{Area}(\Gamma)}{16\pi}}\,  j^{-1/2} + o(j^{-1}) = \frac{\sqrt{\mbox{Area}(\Gamma)}}{4\sqrt{\pi}}\, j^{-1/2} + o(j^{-1})
\eeq
as desired. 
\end{pf}

\subsubsection{Examples of eigenvalues of the single-layer potential operator in three dimensions}\label{sec: Examples of eigenvalues of single-layer potentials in three dimensions}
We present two examples, a sphere and a spherical annulus, for which the eigenvalues can be computed explicitly.
Then we find that a similar phenomenon occurs in the three dimensional single-layer potential operator. 
The explicit eigenvalue formulas for the single-layer potential operator on the sphere reveal nontrivial oscillations 
in the remainder of order $j^{-1}$.  In particular, the remainder estimate $O(j^{-1})$ cannot be improved to $o(j^{-1})$ in this case.

Let $\GO=B(0,R)$ be a ball of radius $R$, and let $\GG=\p\GO$. Let $\{Y_n^m:-n\leq m\leq n,\ n=0,1,2,\ldots\}$ be an orthonormal basis of spherical harmonics in $L^2(\Sbb^2)$, and define the density $f_n^m(y)=Y_n^m(y/R)$ on $\GG$.
Write $x=r\omega$, where $r=|x|$ and $\omega\in\Sbb^2$.
The following single-layer formula is well known (see, e.g., \cite{ACKLM2}):
\begin{align}
\label{eq:single_layer_sphere_action}
\Scal_{\GO}[f_n^m](r\omega)
&=-\frac{R}{2n+1}
\begin{cases}
\displaystyle \left(\frac{r}{R}\right)^nY_n^m(\omega), & r\leq R, \\
\displaystyle \left(\frac{R}{r}\right)^{n+1}Y_n^m(\omega), & r>R,
\end{cases}
\qquad n=0,1,2,\ldots.
\end{align}
Since $\{f_n^m:n\in\Nbb_{\ge0},\ -n\leq m\leq n\}$ is a complete orthogonal basis of $L^2(\GG)$, the above formulas give all eigenvalues of $\Scal_\GG$.
Therefore, for $n=0,1,2,\ldots$,
\begin{align}
\label{eq:sphere_eigenvalue_blocks}
\Gl_j(\Scal_\GG)=-\frac{R}{2n+1},\qquad n^2+1\leq j\leq(n+1)^2.
\end{align}
At the first and last indices of each block, we obtain the following expansions:
\begin{enumerate}
\item[\textup{(1)}] For $j=n^2+1$ ($n=\sqrt{j-1}$), we have
\begin{align}
\label{eq:sphere_first_index_asymptotic}
\Gl_j(\Scal_\GG)
&=-\frac{R}{2\sqrt{j-1}+1}\notag\\
&=-\frac{R}{2\sqrt{j}}+\frac{R}{4j}+O(j^{-3/2}).
\end{align}
\item[\textup{(2)}] For $j=(n+1)^2$ ($n=\sqrt{j}-1$), we have
\begin{align}
\label{eq:sphere_last_index_asymptotic}
\Gl_j(\Scal_\GG)
&=-\frac{R}{2(\sqrt{j}-1)+1}\notag\\
&=-\frac{R}{2\sqrt{j}}-\frac{R}{4j}+O(j^{-3/2}).
\end{align}
\end{enumerate}
In facts, we have 
\begin{align}
- \frac{R}{4} = \liminf_{j \to \infty} j \left(\lambda_j(\Scal_{\Gamma}) + \frac{R}{2\sqrt{j}} \right) <\limsup_{j \to \infty} j \left(\lambda_j(\Scal_{\Gamma}) + \frac{R}{2\sqrt{j}} \right) =\frac{R}{4}
\end{align}

Since every geodesic on a sphere is closed, the nonperiodicity assumption in Theorem \ref{thm: Weyl asymptotic of single-layer potentials in three dimensions} is not satisfied on the sphere $\Gamma={\mathbb S}^2$. Moreover, the explicit formulas of  eigenvalues show that the remainder term of order $j^{-1}$ exhibits oscillatory behavior. Consequently, the improved remainder estimate $o(j^{-1})$ does not hold for the sphere.

Now consider the spherical annulus $\GO = B(0,R_2) \setminus \overline{B(0,R_1)}$, where $0<R_1<R_2$, and set $\Gr:=R_1/R_2$ ($0<\Gr<1$). Its boundary $\GG=\p\GO$ consists of two connected components, the inner sphere $\GG_1$ and the outer sphere $\GG_2$.
We denote the angular densities on each boundary sphere by $f_n^m(y)=Y_n^m(y/|y|)$, $-n\leq m\leq n$, $n\in\Nbb_{\ge0}$.
As in \eqref{eq:boundary_density_identification}, we identify a density $\Gvf$ on $\GG$ with the column vector $[\Gvf_1,\Gvf_2]^T$, where $\Gvf_1:=\Gvf|_{\GG_1}$ and $\Gvf_2:=\Gvf|_{\GG_2}$.
It then follows from \eqref{eq:single_layer_sphere_action} that, for $-n\leq m\leq n$ and $n\in\Nbb_{\ge0}$,
\begin{align}
\Scal_{\GG}\left[\begin{matrix}f_n^m\\0\end{matrix}\right]
&=-\frac{1}{2n+1}\begin{bmatrix}R_1Y_n^m(\omega)\\[4pt]R_1\Gr^{n+1}Y_n^m(\omega)\end{bmatrix}, \\
\Scal_{\GG}\left[\begin{matrix}0\\f_n^m\end{matrix}\right]
&=-\frac{1}{2n+1}\begin{bmatrix}R_2\Gr^nY_n^m(\omega)\\[4pt]R_2Y_n^m(\omega)\end{bmatrix}.
\end{align}
Therefore, for each $n\in\Nbb_{\ge0}$ and $-n\leq m\leq n$, the subspace spanned by $[f_n^m,0]^T$ and $[0,f_n^m]^T$ is invariant under $\Scal_{\GG}$.
With respect to this ordered basis, the restrictions of $\Scal_{\GG}$ are represented by
\begin{align}
\label{eq:spherical_annulus_matrix}
B_n=-\frac{1}{2n+1}
\begin{bmatrix}
R_1 & R_2\Gr^n \\
R_1\Gr^{n+1} & R_2
\end{bmatrix}, \qquad n\in\Nbb_{\ge0},
\end{align}
whose eigenvalues are given by
\begin{align}
\label{eq:spherical_annulus_eigenvalues}
-\frac{R_1+R_2\pm\sqrt{(R_2-R_1)^2+4R_1R_2\Gr^{2n+1}}}{2(2n+1)}.
\end{align}
Since $B_n$ is independent of $m$, each of its eigenvalues contributes multiplicity $2n+1$, corresponding to $m=-n,\ldots,n$.
For $i=1,2$ and $n\in\Nbb_{\ge0}$, we therefore define
\begin{align}
\label{eq:spherical_annulus_component_families}
\Gm_j^{(i)}=-\frac{R_1+R_2+(-1)^i\sqrt{(R_2-R_1)^2+4R_1R_2\Gr^{2n+1}}}{2(2n+1)},
\end{align}
for $n^2+1\leq j\leq(n+1)^2$.
Since $\{[f_n^m,0]^T,[0,f_n^m]^T:n\in\Nbb_{\ge0},\ -n\leq m\leq n\}$ is a complete orthogonal basis of $L^2(\GG)$, the above formulas give all eigenvalues of $\Scal_\GG$.

For $i=1,2$, it follows from \eqref{eq:spherical_annulus_component_families} that
\begin{align}
\Gm_j^{(i)}=-\frac{R_i}{2n+1}+O\left(\frac{\Gr^{2n+1}}{n+1}\right),\qquad n^2+1\leq j\leq(n+1)^2.
\end{align}
At the first and last indices of each block, we obtain the following expansions:
\begin{enumerate}
\item[\textup{(1)}] For $j=n^2+1$ ($n=\sqrt{j-1}$), we have
\begin{align}
\Gm_j^{(i)}
&=-\frac{R_i}{2\sqrt{j-1}+1}
+O\left(\frac{\Gr^{2\sqrt{j-1}+1}}{\sqrt{j}}\right)\notag\\
&=-\frac{R_i}{2\sqrt{j}}+\frac{R_i}{4j}+O(j^{-3/2}).
\end{align}
\item[\textup{(2)}] For $j=(n+1)^2$ ($n=\sqrt{j}-1$), we have
\begin{align}
\Gm_j^{(i)}
&=-\frac{R_i}{2\sqrt{j}-1}
+O\left(\frac{\Gr^{2\sqrt{j}-1}}{\sqrt{j}}\right)\notag\\
&=-\frac{R_i}{2\sqrt{j}}-\frac{R_i}{4j}+O(j^{-3/2}).
\end{align}
\end{enumerate}
Since $0<\Gr<1$, both exponential remainder terms are absorbed into $O(j^{-3/2})$.
Note that
\begin{align}
\Gm_j^{(i)}=-\frac{R_i}{2\sqrt{j}}+O(j^{-1})
=-\frac{\sqrt{\operatorname{Area}(\GG_i)}}{4\sqrt{\pi j}}+O(j^{-1}),\qquad i=1,2.
\end{align}
Although an explicit formula for $\Gl_j(\Scal_{\GG})$ is complicated, we find that
\begin{align}
\Gl_j(\Scal_{\GG})
&=-\frac{\sqrt{R_1^2+R_2^2}}{2\sqrt{j}}+o(j^{-1/2}) \\
&=-\frac{\sqrt{\operatorname{Area}(\GG)}}{4\sqrt{\pi}\sqrt{j}}+o(j^{-1/2}),
\end{align}
since $\operatorname{Area}(\GG)=4\pi(R_1^2+R_2^2)$.
These asymptotics agree with Theorem \ref{thm: Weyl asymptotic of single-layer potentials in three dimensions}.

\subsection{Eigenvalues of the NP operator in three dimensions}\label{sec: Eigenvalues of the NP operator in three dimensions}
We emphasize that the proof here 
is essentially the same as that in Section \ref{sec:  single-layer potentials in three dimensions}. 
\begin{pf}[Proof of Theorem \ref{thm: Weyl's law for NP}] 
Let
\beq 
P:=\widetilde{\mathcal K}_\Gamma =(-\Scal_\Gamma)^{-1/2} \Kcal_\Gamma (-\mathcal S_\Gamma)^{1/2}.
\eeq
Then $P$ is considered as a compact self-adjoint classical $\Gy$DO of order $-1$
on $L^2(\Gamma)$, with principal symbol $K_0(x,\xi)$.

It follows from the assumption in Theorem \ref{thm: Weyl's law for NP} that the second fundamental form of $\Gamma$ is definite at every point, or equivalently, that the Gaussian curvature is everywhere positive. Thus 
$
|K_0(x,\xi)|>0 \ \text{for all}\ (x,\xi)\in T^*\Gamma\setminus 0.
$
We can also decompose $\Gamma$ into $\Gamma^+$ and $\Gamma^-$, on which 
the mean curvature is positive and negative, respectively. We note that $\Gamma^+$ and $\Gamma^-$ 
consist of distinct connected components of $\Gamma$.
Let $P_{\Gamma^+}$ and $P_{\Gamma^-}$ denote the diagonal blocks of $P$ acting 
on $L^2(\Gamma^+)$ and $L^2(\Gamma^-)$, respectively.
Since $\Gamma^+$ and $\Gamma^-$ are disjoint unions of connected components of $\Gamma$, we have
the orthogonal decomposition $L^2(\Gamma)=L^2(\Gamma^+)\oplus L^2(\Gamma^-)$. 
With respect to this decomposition, we can write 
\beq
P=P_{\Gamma^+} \oplus P_{\Gamma^{-}} +R,
\quad R\in\Psi^{-\infty}(\Gamma),
\eeq
where $P_{\Gamma^\pm}$ denote the compressions of $P$ to $L^2(\Gamma^\pm)$, respectively.
We emphasize that the smoothing remainder $R\in\Psi^{-\infty}(\Gamma)$ does not affect
the Weyl asymptotics (see, e.g., \cite[Lemma 2.1]{GPPS}). 
Consequently, it suffices to study the Weyl asymptotics of $P_{\Gamma^+}$ and $P_{\Gamma^-}$ separately.

The operators $P_{\Gamma^\pm}$ are self-adjoint elliptic $\Gy$DOs of order $-1$.
By the existence of elliptic parametrices and the compactness of $\Gamma^\pm$, we obtain
\beq
\dim\ker P_{\Gamma^{\pm}}<\infty.
\eeq
Let $\Pi_{\leq 0}$ and $\Pi_{\geq 0}$ denote the orthogonal spectral projections of $P_{\Gamma^+}$
and $P_{\Gamma^-}$ corresponding to the intervals $(-\infty,0]$ and $[0,\infty)$, respectively.

For $c>0$, we define
\beq
P_{+}:=P_{\Gamma^{-}}(I-\Pi_{\leq 0})+c\Pi_{\leq 0} \ \text{and}\ P_{-}:=P_{\Gamma^{+}}(I-\Pi_{\geq 0})- c\Pi_{\geq 0} .
\eeq
Since $\Pi_{\leq 0}$ and  $\Pi_{\geq 0}$ are finite-rank smoothing operators, we have
\beq
P_{\pm} - P_{\Gamma^{\mp}}   \in\Psi^{-\infty}(\Gamma).
\eeq
Moreover, $P_+$ is positive, self-adjoint, and elliptic of order $-1$ and $P_-$ is negative, self-adjoint, and elliptic of order $-1$. 
Thus we can consider the inverse
\beq
A_{\pm}:=P_{\pm}^{-1}, \qquad D(A_{\pm})=H^1(\Gamma^{\mp}).
\eeq
Then $A_{\pm}$ is a self-adjoint elliptic $\Gy$DO of order $1$ with compact resolvent, and 
\beq
\sigma_{\mathrm{prin}}(A_{\pm})(x,\xi)
=K_0(x,\xi)^{-1} \ \text{on}\ T^{*}\Gamma^{\mp} \setminus \{\xi=0\}.
\eeq
It follows from the vanishing of the subprincipal symbol established in Section \ref{sec: symbol} that
\beq
\operatorname{sub}(A_{\pm})=0.
\eeq
Under the nonperiodicity assumption of Theorem \ref{thm: Weyl's law for NP}, Proposition \ref{thm: Weyl law for psiDO of order 1}  yields
\beq
N(A_{\pm},\lambda)=\sharp \{j\in \Nbb\ |\  |\lambda_{j}(A_{\pm})| < \lambda\} =  \alpha_{\pm}\lambda^2+o(\lambda),
\qquad \lambda\to\infty,
\eeq
where
\beq
\alpha_{\pm} =\frac{1}{(2\pi)^2} \operatorname{Vol} \left\{ (x,\xi)\in T^*\Gamma^{\mp}: |K_0(x,\xi)^{-1}|<1 \right\} 
= \frac{3W(\Gamma^{\mp})-2\pi\chi(\Gamma^{\mp})}{128\pi} 
\eeq
Here $W(\Gamma)$ and $\chi(\Gamma)$ denote the Willmore energy and Euler characteristic of $\Gamma$, respectively
(see, e.g., \cite{Miyanishi:Weyl}). 

Taking account of the negative and positive eigenvalues, the refined Weyl's law implies
\beq
\lambda_j(A_{\pm}) = \pm \sqrt{\frac{j}{\alpha_{\pm}}}+o(1).
\eeq

Taking reciprocals, we obtain
\beq
\lambda_j(P_{\pm}) = \pm \sqrt{\alpha_\pm}\,j^{-1/2} +o(j^{-1}).
\eeq

Since $P_{\pm}$ differs from $P_{\Gamma^{\mp}}$ only on a finite-dimensional spectral subspace,
this asymptotic formula also holds for the eigenvalues of $P$.
In fact, a fixed shift in the eigenvalue indices changes the leading term only by $O(j^{-3/2})$, which is $o(j^{-1})$.

Finally, $P$ and $\mathcal K_\Gamma$ have the same
nonzero eigenvalues, counted with multiplicity.

Since $P_+$ and $P_-$ have only finitely many
nonpositive and nonnegative eigenvalues, respectively,
the asymptotic behavior of the positive and negative
eigenvalues is determined by $P_+$ and $P_-$,
respectively.

Consequently, under the nonperiodicity assumption, we obtain
\beq
\lambda_j^{+}(\mathcal K_\Gamma) =
\left(\frac{3W(\Gamma^{-})-2\pi\chi(\Gamma^{-})}
{128\pi}\right)^{1/2}j^{-1/2}+o(j^{-1}).
\eeq
Similarly,
\beq
\lambda_j^{-}(\mathcal K_\Gamma)= -\left(
\frac{3W(\Gamma^{+})-2\pi\chi(\Gamma^{+})}{128\pi} \right)^{1/2} j^{-1/2}+o(j^{-1}).
\eeq
This completes the proof.
\end{pf}
\subsubsection{Examples of eigenvalues of the NP operator in three dimensions}\label{sec: Examples of eigenvalues of NP operator in three dimensions}
As is known in \cite{Miyanishi:Weyl, MS} and references therein, the explicit formula for eigenvalues of the NP operator on the unit sphere reads 
\[
\Gl_j(\Kcal_\GG) = \frac{1}{2(2n+1)}, \quad n^2 + 1 \leq j \leq (n+1)^2, \quad n = 0, 1, 2, \ldots.
\]
By the expansion of the integral kernel of the NP operator with respect to the spherical harmonics, it can be seen that the corresponding eigenfunctions are the spherical harmonics. In the same way as for the asymptotic behaviors \eqref{eq:sphere_first_index_asymptotic} and \eqref{eq:sphere_last_index_asymptotic} of the eigenvalues of the single-layer potential operator, we obtain
\begin{enumerate}
\item[\textup{(1)}] For $j=n^2+1$ ($n=\sqrt{j-1}$), we have
\begin{align*}
\Gl_j(\Kcal_\GG) =\frac{1}{4\sqrt{j-1}+2} =\frac{1}{4\sqrt{j}} - \frac{1}{8j} + O(j^{-3/2}).
\end{align*}
\item[\textup{(2)}] For $j=(n+1)^2$ ($n=\sqrt{j}-1$), we have
\begin{align*}
\Gl_j(\Kcal_\GG) =\frac{1}{4(\sqrt{j}-1)+2} =\frac{1}{4\sqrt{j}} + \frac{1}{8j}+O(j^{-3/2}).
\end{align*}
\end{enumerate}
We note that, in this case, the NP flow coincides with the geodesic flow up to a constant rescaling of time, and the set of periodic points has positive Liouville measure. This situation is the same as that for the single-layer potential operator discussed in Section \ref{sec: Examples of eigenvalues of single-layer potentials in three dimensions}. Consequently, the improved remainder estimate in Theorem \ref{thm: Weyl's law for NP} cannot be applied.

As a second example, we consider the spherical shell $\GO=B(0,R_2)\setminus\overline{B(0,R_1)}$, where $0<R_1<R_2$.
We define $\Gr:=R_1/R_2$.  The boundary $\GG=\partial\GO$ consists of two concentric spheres, 
the inner sphere $\GG_1$ and the outer sphere $\GG_2$. 
The eigenvalues of the NP operator on $\GG$ are given by the following formulas
(see, e.g., \cite{ACKLM2}). 
\beq\label{eq:spherical_annulus_NP_eigenvalue_blocks} 
\lambda_j^\pm(\Kcal_\GG) = \pm\frac{\sqrt{1+4n(n+1)\Gr^{2n+1}}}{2(2n+1)},
\quad n^2+1\leq j\leq(n+1)^2,
\quad n=0,1,2,\ldots.
\eeq
with the ordering used in Theorem \ref{thm: Weyl's law for NP}.  Since $0<\Gr<1$, the difference between \eqref{eq:spherical_annulus_NP_eigenvalue_blocks} and $\pm1/[2(2n+1)]$ is $O(j^{-\infty})$. On each boundary component, the NP flow coincides with
the geodesic flow up to a constant rescaling of time. Thus all trajectories are periodic, and the nonperiodicity
assumption in Theorem \ref{thm: Weyl's law for NP} is not satisfied. 
As in the case of the sphere, the remainder of order $j^{-1}$ exhibits nontrivial oscillations, 
and the estimate $O(j^{-1})$ cannot be improved to $o(j^{-1})$ for either eigenvalue sequence.

On the other hand, explicit eigenvalue formulas for the NP operator are known for certain smooth surfaces 
in three dimensions, including ellipsoids (see, e.g., \cite{FengKang2016, Martensen1999, MS, Ritter1995} and references therein).
For non-spherical ellipsoids and generic smooth strictly convex surfaces, we expect the improved asymptotic estimate in Theorem \ref{thm: Weyl's law for NP} to hold. However, we have not yet been able to verify the required nonperiodicity condition for the NP flow in these cases. We leave this problem for future investigation.

{\color{black}
\section{Summary}\label{sec: summary}

In this paper, we have investigated the asymptotic behavior of the eigenvalues of two important classes of layer potential operators: the single-layer potential operator and the Neumann--Poincar\'e (NP) operator in two and three dimensions.

For the single-layer potential operator, the leading asymptotic behavior of the eigenvalues is essentially governed by the size of the boundary. In two dimensions, we obtain precise eigenvalue asymptotics with superpolynomially decaying remainders,
determined by the lengths of the boundary components. In three dimensions, we have obtained a more precise asymptotic formula under the assumption that the set of periodic points of the geodesic flow on the unit cotangent bundle of the boundary has measure zero.

For the NP operator, we have focused primarily on the three-dimensional case. In particular, for smooth surfaces with everywhere positive Gaussian curvature, we have refined the Weyl asymptotics established in previous studies by improving the remainder estimate from $o(j^{-1/2})$ to $o(j^{-1})$, under the assumption that the set of periodic points of the associated Hamiltonian flow, which we refer to as the NP flow, has Liouville measure zero.

These results provide further evidence that the geometry of the boundary plays a fundamental role in governing the spectral asymptotics of layer potential operators. In this sense, the present work can be regarded as a study of the spectral geometry of layer potential operators.}


\end{document}